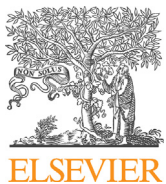
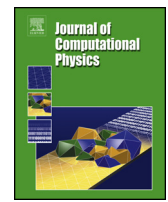

# High-order numerical solutions to the shallow-water equations on the rotated cubed-sphere grid

Stéphane Gaudreault [a,*], Martin Charron [a], Valentin Dallerit [b], Mayya Tokman [b]

[a] *Recherche en prévision numérique atmosphérique, Environnement et Changement climatique Canada, 2121 Route Transcanadienne, Dorval, Québec, H9P 1J3, Canada*
[b] *School of Natural Sciences, University of California, 5200 N. Lake Road, Merced, CA 95343, United States*

**A B S T R A C T**

A novel numerical approach to solving the shallow-water equations on the sphere using high-order numerical discretizations in both space and time is proposed. A space-time tensor formalism is used to express the equations of motion covariantly and to describe the geometry of the rotated cubed-sphere grid. The spatial discretization is done with the direct flux reconstruction method, which is an alternative formulation to the discontinuous Galerkin approach. The equations of motion are solved in differential form and the resulting discretization is free from quadrature rules. It is well known that the time step of traditional explicit methods is limited by the phase velocity of the fastest waves. Exponential integration is employed to enable integrations with significantly larger time step sizes and improve the efficiency of the overall time integration. New multistep-type exponential propagation iterative methods of orders 4, 5 and 6 are constructed and applied to integrate the shallow-water equations in time. These new schemes enable time integration with high-order accuracy but without significant increases in computational time compared to low-order methods. The exponential matrix functions-vector products used in the exponential schemes are approximated using the complex-step approximation of the Jacobian in the Krylov-based KIOPS (Krylov with incomplete orthogonalization procedure solver) algorithm. Performance of the new numerical methods is evaluated using a set of standard benchmark tests.

## 1. Introduction

A challenge of great interest in the numerical weather prediction community is to develop numerical methods for the governing equations that are conservative, highly accurate, geometrically flexible, computationally efficient, and simply formulated. Second-order numerical methods are often preferred in operational models due to their simplicity and robustness. High-order methods are sometimes perceived as less robust and complicated to implement. The aims of this paper are to present: 1) spatial numerical techniques that are almost as simple as finite differences; and 2) multistep exponential time integration methods with superior stability properties compared to explicit schemes.

* Corresponding author.
 *E-mail addresses:* stephane.gaudreault@ec.gc.ca (S. Gaudreault), martin.charron@ec.gc.ca (M. Charron), vdallerit@ucmerced.edu (V. Dallerit), mtokman@ucmerced.edu (M. Tokman).

For spatial discretization aspects, the direct flux reconstruction (DFR) scheme [1] is applied to the cubed-sphere grid. The DFR method is a simplified formulation of a class of methods called flux reconstruction (FR) [2]. The DFR method is equivalent to the weak form of the nodal discontinuous Galerkin (NDG) scheme in one dimension and on multidimensional elements with tensor product bases [1,3]. In this alternative formulation of the NDG method, the conservative equations are solved in differential form and the discretization is free from quadrature rules, resulting in a simple and computationally efficient algorithm. This method is also well suited to large geophysical fluid dynamics problems because it is local and has good conservation properties [4–6].

Time integration with high-order methods is a much more difficult problem. The implementation of high-order Eulerian techniques in operational models is made difficult by the Courant-Friedrichs-Lewy (CFL) condition that limits the time step in several explicit time integration schemes [7]. This motivated the recent investigation of exponential time integrators in geophysical applications [8–17]. These approaches allow for longer time steps and yield higher accuracy than traditional algorithms. In addition, these studies have shown that exponential integrators accurately calculate the entire spectrum of waves propagating in the atmosphere. Beside the obvious advantage of high accuracy, these methods also eliminate the need to divide the governing equations into linear and non-linear parts. This task is performed automatically by evaluating the action of the Jacobian operator. This consideration could greatly simplify the design of numerical models.

In this paper, new high-order exponential propagation iterative (EPI) methods are introduced and tested in combination with DFR for the shallow-water equations on the sphere. The paper is organized as follows. The governing equations in arbitrary coordinates are introduced in §2. This sets the stage for §3 and §4 where the numerical schemes are described. Results from a series of numerical experiments are discussed in §5. Conclusions and future works are presented in §6.

## 2. General form of the equations of motion and the rotated cubed-sphere grid

Solving the shallow-water equations is often the first step towards the construction of a comprehensive numerical weather prediction model. Shallow-water dynamics represents a (2+1)-dimensional (i.e. two spatial and one temporal dimensions) nonlinear system of hyperbolic partial differential equations. They are obtained from the Euler equations by assuming that the depth of the fluid is small compared to the mean radius of the Earth, incompressibility, hydrostasy, weak stratification, and in the case of spherical geometry, that the spherical geopotential approximation holds.

Here, the shallow-water equations will be written in a space-time tensorial formalism. One advantage of using space-time tensor calculus and tensor-related objects is that the equations of motion may be expressed in arbitrary (possibly non-inertial and time-dependent) coordinates within a unified framework.

The synchronously covariant shallow-water equations for continuity and momentum in quasi-flux form are respectively

$$\frac{\partial}{\partial t}\left(\sqrt{g}H\right) + \frac{\partial}{\partial x^j}\left(\sqrt{g}Hu^j\right) = 0, \tag{1}$$

$$\frac{\partial}{\partial t}\left(\sqrt{g}Hu^i\right) + \frac{\partial}{\partial x^j}\left(\sqrt{g}\left[Hu^iu^j + \frac{1}{2}g_rh^{ij}H^2\right]\right) =$$
$$-2\sqrt{g}\,\Gamma^i_{j0}Hu^j - \sqrt{g}\,\Gamma^i_{jk}\left(Hu^ju^k + \frac{1}{2}g_rh^{jk}H^2\right) - \sqrt{g}Hg_rh^{ij}\frac{\partial h_B}{\partial x^j}, \tag{2}$$

where $H$ is the fluid's thickness scalar field, $h_B$ the height of the bottom orography, $u^i$ ($i=1,2$) the two components of the velocity field, $g_r$ the constant effective gravitational acceleration, $g$ the determinant of the covariant space-time metric tensor, $h^{ij}$ the space-only contravariant metric tensor, and $\Gamma$'s the Christoffel symbols in space-time. Notice that the Coriolis effect is associated with the Christoffel symbols $\Gamma^i_{j0}$. The indices $i,j$ take the values 1 and 2, and repeated indices are summed. The derivation of these equations is presented in Appendix A. The form of the governing equations provided by Eqs. (1)–(2) will be discretized in the following sections.

The relevant space-time metric terms associated with the cubed-sphere coordinates on a rotating 2-sphere with radius $a$ and constant angular velocity $\Omega$ are provided in Appendix B. It is assumed that the reader is familiar with the characteristics of the cubed-sphere grid, otherwise [18–28] provide descriptions and details on its properties. In practical applications, it may be desirable to rotate a coordinate system to avoid the co-location of certain geographical points and grid peculiarities such as coordinate discontinuities. In Appendix B, an arbitrary rotation of the cubed-sphere coordinates is considered.

Coordinate lines at the interface of two panels may not be smooth as a result of the composite character of the global cubed-sphere coordinates. Although scalars have by definition the same value at a point on the interface of two panels with different coordinates, higher-rank tensor components are different if expressed with the coordinate basis of one or the other side of the interface. In Appendix C, consistency relations between all interfacing panels are provided for contravariant and covariant first-rank tensors.

## 3. Spatial discretization

The discontinuous Galerkin method is built on a weak integral formulation of the governing equations, where all unknown functions are approximated by high-order polynomials (see [29,30] for a review). This leads to the evaluation of several integrals using quadrature rules. The recent direct flux reconstruction method [1,31] is formulated using a differential form of the equations. Therefore, this method only requires Lagrange polynomials and avoids the need for quadrature rules, which simplifies its implementation. It has been shown that the DFR method results in a scheme equivalent to the weak form of the nodal discontinuous Galerkin method in one dimension and on multidimensional elements with tensor product bases [1,3]. It should be noted that these ideas are not entirely new in geophysical applications since the 3rd order scheme presented in [32] may also be considered as a special case of the DFR method. In §3.1, an overview of the basic properties of the method in one dimension is presented and its extension to two dimensions and curved geometry is discussed in §3.2.

*3.1. Direct flux reconstruction in one dimension*

Consider the following conservation law:

$$\frac{\partial q}{\partial t} + \frac{\partial f(q)}{\partial x} = 0, \quad x \in \mathcal{D}. \tag{3}$$

The calculation domain $\mathcal{D}$ is divided into $N_e$ non-overlapping elements $\mathcal{D}_j = [x_{j-\frac{1}{2}}, x_{j+\frac{1}{2}}]$, for $j = 1, ..., N_e$. The size of the element $j$ is $\Delta_j = x_{j+\frac{1}{2}} - x_{j-\frac{1}{2}}$. In this subsection, all indices refer to positions in a one-dimensional grid and not to space-time indices as in §2. Also, implicit summation over repeated indices is not assumed.

To conveniently treat a general non-uniform grid, a local coordinate variable $\xi$ is defined by an affine transformation mapping each element $\mathcal{D}_j$ onto the so-called reference element $I_j = [-1, 1]$:

$$\xi(x) = \frac{2}{\Delta_j}(x - x_{j-\frac{1}{2}}) - 1 \tag{4}$$

with inverse

$$x(\xi) = \left(\frac{1-\xi}{2}\right) x_{j-\frac{1}{2}} + \left(\frac{1+\xi}{2}\right) x_{j+\frac{1}{2}}. \tag{5}$$

Derivatives in global and local coordinate systems are related by the chain rule:

$$\frac{\partial}{\partial x} = \frac{2}{\Delta_j} \frac{\partial}{\partial \xi}. \tag{6}$$

Similarly to several other high-order methods such as FR and NDG, a set of $N_s \geq 1$ solution points are defined within the reference element. These points are denoted $\xi_k$, for $k = 1, \ldots, N_s$. There are several ways to place these points, leading to schemes with different numerical properties [33]. In the present study, the Gauß-Legendre points have been chosen based on theoretical results which show the optimal nature of these solution points [34,35]. The placement of the Gauß-Legendre points corresponds to the roots of the Legendre polynomials. Although the solution and flux points are collocated, the method does not suffer from the so-called $2\Delta x$ numerical mode that typically affects the unstaggered finite difference schemes on uniform grids. In effect, the entries of the Gauß-Legendre differentiation matrices are all non-zero. This is in contrast with centered unstaggered finite differences on a uniform grid, where the coefficients vanish at the center of the stencil. The Runge phenomenon is also significantly mitigated by this choice of points.

The essence of the DFR method consists in approximating the exact solution $q$ with an approximate solution $\hat{q}$. This approximate solution is obtained by taking the piecewise sum of $N_e$ functions $\hat{q}_j$, each function being defined as a polynomial of degree $N_s - 1$ within the element $\mathcal{D}_j$ and exactly zero elsewhere. Thus, the approximate solution for a given element $j$ can be represented on a reference element as

$$\hat{q}_j(\xi, t) = \sum_{k=1}^{N_s} \hat{q}_{j,k}(t) \ell_k(\xi), \tag{7}$$

where $\hat{q}_{j,k}(t)$ are coefficients defined at the solution points, and

$$\ell_k(\xi) = \prod_{\substack{l=1 \\ l \neq k}}^{N_s} \frac{\xi - \xi_l}{\xi_k - \xi_l} \tag{8}$$

are the Lagrange basis polynomials.

Similarly, the exact flux $f(q)$ in Eq. (3) is approximated by the piecewise sum of polynomials $\hat{f}_j$ of degree $N_s - 1$ within each element and zero outside:

$$\hat{f}_j(\xi, t) = \sum_{k=1}^{N_s} \hat{f}_{j,k}(t)\, \ell_k(\xi), \tag{9}$$

where $\hat{f}_{j,k}(t) = f(\hat{q}_{j,k}(t))$. It is worth noting that even if $\hat{q}$ could be represented exactly with the basis functions, the term $\hat{f}$ is obtained from products of dynamical variables and the expansion (9) would not be exact in general. This could introduce aliasing errors which may require filtering.

At a given time $t$, the functions $\hat{q}$ and $\hat{f}$ are continuous within an element but are not uniquely defined at the element boundaries (i.e., at the points $\xi = \pm 1$). If for instance the flux at $\xi = -1$ on the element $j$ is considered, its value may differ from the value at $\xi = 1$ on the element $j-1$. These disparities constitute a Riemann problem. For smooth solutions, a simple upwinding might be an acceptable approximation for the common flux at the interface. If, however, the physical solution should include sharp gradients, the common flux is generally prescribed by an approximate Riemann solver borrowed from the finite volume methodology (see e.g. [36,37]). The study [25] adapted the Rusanov [38], Roe [39] and AUSM+-up [40] methods for the shallow-water equations. Their simulations showed that AUSM+-up provides the best overall accuracy when applied to various test cases, followed closely by the Roe solver. The Rusanov solver showed significantly worse performance in terms of accuracy and conservation error. However, the AUSM+-up method is much more computationally intensive than the Rusanov method, which explains why the latter is a more popular choice with discontinuous Galerkin methods. After exploring different possibilities, an adaptation of the basic AUSM solver [41] was found to be a good compromise between simplicity, accuracy and numerical efficiency. This Riemann solver is described in Appendix D.

To account for the interaction between adjacent elements, a continuous polynomial, denoted $F_j(\xi)$, is defined within each element. This continuous flux function must satisfy the following conditions:

1. $F_j(\xi)$ must take the value of the Riemann fluxes at both ends of the element $\mathcal{D}_j$:

$$F_j(-1) = \bar{f}\left(q_{j-\frac{1}{2}}^L, q_{j-\frac{1}{2}}^R\right), \quad F_j(1) = \bar{f}\left(q_{j+\frac{1}{2}}^L, q_{j+\frac{1}{2}}^R\right); \tag{10}$$

2. $F_j(\xi)$ must be a polynomial of degree $N_s + 1$;
3. $F_j(\xi)$ must approximate the discontinuous flux function $\hat{f}_j(\xi, t)$ at a given time $t$.

Condition 1 is necessary for the polynomial to be $C^0$. Condition 2 implies that the derivative of $F_j(\xi)$ is a polynomial of degree $N_s$ and that the method has an accuracy of order $N_s$. To satisfy these two conditions, one defines the continuous polynomial over the extended set of solution points as the union of the set of interior solution points and the interface flux points, that is $\{\xi_0, \xi_1, \ldots, \xi_{N_s}, \xi_{N_s+1}\}$, with $\xi_0 = -1$, $\xi_{N_s+1} = 1$ and the usual Gauß-Legendre points for indices $1, \ldots, N_s$. Condition 3 may be satisfied by imposing that $F_j(\xi_k) = \hat{f}_{j,k}$ on the interior Gauß-Legendre points $k = 1, \ldots, N_s$. These three conditions are sufficient to define the following continuous flux polynomial:

$$F_j(\xi) = F_j(-1)\tilde{\ell}_0(\xi) + \left[\sum_{k=1}^{N_s} \hat{f}_{j,k}(t)\, \tilde{\ell}_k(\xi)\right] + F_j(1)\tilde{\ell}_{N_s+1}(\xi), \tag{11}$$

where $\tilde{\ell}_n$ are the Lagrange interpolation polynomial basis constructed from Eq. (8) but now for the $N_s + 2$ interpolation points.

Finally, the derivative of the continuous flux function is given by

$$F_j'(\xi) = F_j(-1)\tilde{\ell}_0'(\xi) + \left[\sum_{k=1}^{N_s} \hat{f}_{j,k}(t)\, \tilde{\ell}_k'(\xi)\right] + F_j(1)\tilde{\ell}_{N_s+1}'(\xi), \tag{12}$$

with

$$\tilde{\ell}_n'(\xi) = \sum_{\substack{i=0 \\ i \neq n}}^{N_s+1} \left[ \frac{1}{\xi_n - \xi_i} \prod_{\substack{m=0 \\ m \neq (i,n)}}^{N_s+1} \frac{\xi - \xi_m}{\xi_n - \xi_m} \right]. \tag{13}$$

After replacing variables and the spatial derivative in Eq. (3) by their discrete counterparts, the following semi-discrete ordinary differential equation is obtained:

$$\frac{d}{dt}\hat{q}_{j,k} = -\frac{2}{\Delta_j} F_j'(\xi_k). \tag{14}$$

In practice, the derivative of the continuous flux function is calculated as a matrix-vector product using a polynomial derivative matrix [42]. The equivalence of this scheme with the weak form of the NDG method may be verified easily by noting that: 1) $\tilde{\ell}'_0(\xi)$ and $\tilde{\ell}'_k(\xi)$ respectively take the same values as left and right Radau polynomials at the Gauß-Legendre points [3], and 2) the part of the differentiation matrix applied to the solution points is identical to, say, Eq. (43) in [43]. Given this equivalence between DFR and NDG, one expects that the numerical properties of the latter, which are well studied in the literature (see e.g. [29,30]), also apply to the former.

### 3.2. Extension to two dimensions and curved geometry

The extension of the DFR method to quadrilateral and hexahedral elements, for which the basis polynomials are the tensor product of the one-dimensional basis functions, is straightforward. The basic idea is to first transform the governing equations from a physical element to the reference or standard element. Then, the one-dimensional DFR formulation is applied in each coordinate direction of the standard element without further complications. The discretization of the tensorial momentum equations is performed by treating each component separately.

Although the method is remarkably similar in flat and curved geometries, it is possible to simplify the formulation by avoiding the transformation of the derivatives given in Eqs. (6) and (14). For instance, in Eqs. (1) and (2) one may take advantage of the tensorial formulation by transforming the velocities $u^i$, metric $h^{ij}$ and factor $\sqrt{g}$ as well as the components of the Christoffel symbols $\Gamma^i_{0k}$ and $\Gamma^i_{jk}$ from the cubed-sphere coordinates to the local coordinate system $[-1,1] \times [-1,1]$ on each element. The transformation rules are

$$\tilde{u}^j = \frac{2}{\Delta x^j} u^j, \tag{15}$$

$$\sqrt{\tilde{g}} = \frac{1}{4} \Delta x^1 \Delta x^2 \sqrt{g}, \tag{16}$$

$$\tilde{g}_{ij} = \frac{1}{4} \Delta x^i \Delta x^j g_{ij}, \tag{17}$$

$$\tilde{h}^{ij} = \frac{4}{\Delta x^i \Delta x^j} h^{ij}, \tag{18}$$

$$\tilde{\Gamma}^i_{jk} = \frac{1}{2} \Delta x^i \Gamma^i_{jk}, \tag{19}$$

$$\tilde{\Gamma}^i_{0k} = \Gamma^i_{0k}, \tag{20}$$

where $\Delta x^1$ and $\Delta x^2$ are the element dimensions in the cubed-sphere coordinate system and the tilde symbol indicates the corresponding values in the coordinate system of the standard element. These transformations may be computed only once at the beginning of the simulation and then used consistently in all calculations. Note that these transformations do not affect the consistency relations at the interfaces of the cubed-sphere panels.

## 4. Time integration

The semi-discrete system obtained after applying the DFR to the shallow-water equations is an example of a general problem arising in the so-called method of lines, in which the space derivatives of a partial differential equation (PDE) are discretized first, leading to a large autonomous system of ordinary differential equations (ODEs) of the form

$$\frac{d}{dt} q(t) = \mathcal{F}(q(t)), \quad q(t_0) = q_0, \tag{21}$$

where $q(t)$ represents the unknown dynamical quantities and $\mathcal{F}$ is a function describing all forcing terms driving the system.

Numerical time integration of Eq. (21) is a challenging task due to both high dimensionality of $\mathcal{F}(q)$ as well as the stiffness of its Jacobian $\partial \mathcal{F}/\partial q$. The latter property is due to the presence of a wide range of temporal frequencies characteristic of the evolution of this system. Traditionally, there have been two general approaches to temporal discretization of Eq. (21): explicit and implicit methods (see e.g. [44]). On the one hand, explicit methods are very efficient per time step since such algorithms require only evaluations of the right-hand side function $\mathcal{F}(q)$ using precomputed values of $q$ and possible additions of such vectors. However, explicit integrators suffer from poor numerical stability properties, and for stiff systems the stability constraints impose a time step that is often too small for practical applications. Implicit methods, on the other hand, usually possess much better stability properties and allow time integration of the system using significantly larger time step sizes compared to explicit schemes. Such improvement though is accompanied by an increase in computational complexity per time step since implicit methods require solution of a large and stiff system of nonlinear equations at each time step. Newton-Krylov algorithms [45] are the most commonly used in practice to solve such nonlinear systems with implicit time stepping. The nonlinearity and stiffness of many dynamical systems, such as the shallow-water equations, may result in slow convergence of both the Krylov method used to solve the underlying linear system within a Newton iteration and the Newton iteration itself. Preconditioning is used to alleviate this challenge, but developing an effective preconditioner for many problems can in itself be a difficult and time-consuming task. This is particularly true when higher-order spatial discretizations, such as the DFR method, are employed to approximate spatial operators. Given these considerations, there is a pressing need for new time integration techniques with better stability properties compared to explicit methods and improved computational complexity per time step compared to implicit schemes. In recent years, exponential integration emerged as a possible alternative to implicit methods in solving large stiff systems.

Let $\{t_n\}_{n=0}^{M}$ be a set of nodes that represent some discretization of the integration interval $[t_0, t_{end}]$. The numerical solution of Eq. (21) at time $t_n$ is denoted as $q_n \approx q(t_n)$ and the next task is to develop a time integration method to compute the vectors $q_n$. The starting point for constructing an exponential integrator is often a Taylor expansion of $\mathcal{F}(q)$ around a known value of the solution $q_n$ that allows writing the equation in the following form

$$\frac{dq}{dt} = \mathcal{F}(q) = \mathcal{F}(q_n) + \frac{\partial \mathcal{F}}{\partial q}(q_n)(q - q_n) + \mathcal{R}(q), \tag{22}$$

where $\mathcal{R}(q) = \mathcal{F}(q) - \mathcal{F}(q_n) - \frac{\partial \mathcal{F}}{\partial q}(q_n)(q - q_n)$ is the nonlinear remainder after the first two terms of the Taylor expansion. Denoting the Jacobian matrix $\mathcal{J}_n = \frac{\partial \mathcal{F}}{\partial q}(q_n)$ for brevity and using an integrating factor $e^{t\mathcal{J}_n}$, one may write the integral form of Eq. (22), and the solution at $t_n + \Delta t$ is given by the variation of constants formula

$$q(t_n + \Delta t) = q_n + \mathcal{J}_n^{-1}(e^{\Delta t \mathcal{J}_n} - I)\mathcal{F}(q_n) + \int_{t_n}^{t_n + \Delta t} e^{(t_n + \Delta t - t)\mathcal{J}_n} \mathcal{R}(q(t))dt \tag{23}$$

or, if one sets $\varphi_1(z) = (e^z - 1)/z$, as

$$q(t_n + \Delta t) = q_n + \varphi_1(\Delta t \mathcal{J}_n)\Delta t \mathcal{F}(q_n) + \int_{t_n}^{t_n + \Delta t} e^{(t_n + \Delta t - t)\mathcal{J}_n} \mathcal{R}(q(t))dt. \tag{24}$$

An exponential integrator to estimate $q(t_n + \Delta t)$ is then constructed by approximating the nonlinear integral in Eq. (24) to the desired order of accuracy. If the integral is neglected, one obtains a second-order method known as an Exponential Euler scheme (hereafter denoted EPI2 as in [46]):

$$q(t_n + \Delta t) = q_n + \varphi_1(\Delta t \mathcal{J}_n)\Delta t \mathcal{F}(q_n). \tag{25}$$

To increase the order of accuracy, the nonlinear integral is usually approximated using Runge-Kutta or multistep-type approaches. Either of such schemes will result in using some form of a polynomial approximation to the function $\mathcal{R}(q(t))$ in variable $t$ to estimate the nonlinear integral in Eq. (24). As a result, the approximate solution $q_{n+1} \approx q(t_n + \Delta t)$ will be expressed in terms of products of matrix functions and vectors like $\varphi_k(\Delta t \mathcal{J}_n)\mathcal{R}_k$, where $\mathcal{R}_k$ is a vector obtained by computing the nonlinear function $\mathcal{R}(q_k)$ using a known value of $q_k$. The functions $\varphi_k$ satisfy the relation

$$\varphi_0(z) = e^z, \quad \varphi_k(z) = \int_0^1 e^{(1-s)z} \frac{s^{k-1}}{(k-1)!}ds, \quad k \geq 1. \tag{26}$$

All $\varphi_k(z)$ are analytic functions defined on the complex plane. This definition can then be extended to matrices using any of the available definitions of matrix functions (see e.g. [47]). For instance, they could be defined using power series similar to the matrix exponential:

$$\varphi_k(A) = \sum_{n=0}^{\infty} \frac{A^n}{(n+k)!}.$$

Since $\mathcal{J}_n$ is generally a large stiff matrix, evaluation of the products $\varphi_k(\Delta t \mathcal{J}_n)\mathcal{R}_k$ is the most computationally intensive part of approximating $q_{n+1}$ at every time step. Thus, choosing the right algorithm to approximate these products is key to implementing an efficient exponential integrator.

While there are a number of algorithms that have been proposed for estimating products of matrix functions and vectors, only a few of those are suitable for general large stiff problems, and even fewer are appropriate for problems where the argument matrix is not known or cannot be stored explicitly. Many of such algorithms are only computationally feasible for small matrices [48,49]. Other methods can be used for large matrices but require some information about the norm or the spectrum of the matrices [50–52]. For many systems including the one considered here, it is, however, common to have the Jacobian matrix only available in the so-called matrix-free form—i.e. only a function that multiplies the Jacobian matrix with a vector is available rather than having an explicitly stored Jacobian. When an approximation to the Jacobian-vector product is used, the resulting method is then called an inexact Krylov method. An analysis of inexact Krylov subspace methods for approximating the matrix exponential has been presented in [53].

A popular matrix-free method (see [45] and references therein) is the finite-difference approximation

$$\mathcal{J}_n v = \frac{\mathcal{F}(q_n + \epsilon \cdot v) - \mathcal{F}(q_n)}{\epsilon} + \mathcal{O}(\epsilon), \tag{27}$$

where $\epsilon \in \mathbb{R}$ is a small parameter and $v$ is an arbitrary vector. It is well known that if $\epsilon$ is too small, then the rounding errors incurred in the evaluation of $\mathcal{F}$ begin to dominate. For this reason, a methodology based on the complex-step approximation is rather used [54]:

$$\mathcal{J}_n v = \Im\left[\frac{\mathcal{F}(q_n + \epsilon\, i \cdot v)}{\epsilon}\right] + \mathcal{O}(\epsilon^2), \tag{28}$$

where $i^2 = -1$ and $\Im$ denotes the imaginary part operator.

Hence, the result of the Jacobian-vector product is calculated using only one evaluation of the right-hand side function with complex arguments. The complex-step approximation has gained popularity in the field of machine learning in recent years [55] but has seldom been used in the context of exponential integrators.

For large stiff matrix-free cases, Krylov-projection-type algorithms are the best approach to estimating products $\varphi_k(\Delta t \mathcal{J}_n) \mathcal{R}_k$. In particular, adaptive Krylov-based methods [56], including the recently proposed KIOPS algorithm [10], have been shown to deliver the most efficiency for such problems. Given a matrix $A$ and a set of vectors $b_0, ..., b_p$ these methods allow fast computation of the linear combinations of the type

$$\varphi_0(A)b_0 + \varphi_1(A)b_1 + \varphi_2(A)b_2 + ... + \varphi_p(A)b_p. \tag{29}$$

The adaptive Krylov-based methods compute the expression (29) using only one Krylov projection and an adaptive substepping mechanism which speeds up the computation by appropriately scaling the matrix $A$. The precise formulation of the KIOPS algorithm (its structure and implementation) is given in [10].

As stated earlier, obtaining the matrix functions-vector products $\varphi_k(\Delta t \mathcal{J}_n) \mathcal{R}_k$ and their linear combinations represents the main computational load of an exponential integrator. One should therefore minimize the number of such operations and optimize their computation. Exponential Propagation Iterative (EPI) methods of Runge-Kutta [57] and multistep [46] types were designed to achieve this. Simulations with weather and climate models usually employ constant time steps because they are controlled not only by the resolved dynamics but also by constraints on sub-grid scale parameterizations. Therefore, only methods with constant time steps are considered here.

A third-order multistep-type EPI method is derived in [46]. The EPI3 algorithm to solve Eq. (22) may be written as

$$q_{n+1} = q_n + \varphi_1(\Delta t \mathcal{J}_n)\Delta t \mathcal{F}_n + \tfrac{2}{3}\varphi_2(\Delta t \mathcal{J}_n)\Delta t \mathcal{R}_{n-1}, \tag{30}$$

where $\mathcal{F}_n = \mathcal{F}(q_n)$ and

$$\varphi_1(z) = \frac{e^z - 1}{z}, \tag{31}$$

$$\varphi_2(z) = \frac{e^z - 1 - z}{z^2}, \tag{32}$$

$$\mathcal{R}_{n-1} = \mathcal{F}(q_{n-1}) - \mathcal{F}_n - \mathcal{J}_n(q_{n-1} - q_n). \tag{33}$$

In the following, higher-order multistep-type EPI algorithms are presented. The following general ansatz will become useful:

$$q_{n+1} = q_n + \varphi_1(\Delta t \mathcal{J}_n)\Delta t \mathcal{F}(q_n) + \sum_{m=1}^{M} \varphi_m(\Delta t \mathcal{J}_n) v_m, \tag{34}$$

$$v_m = \sum_{i=1}^{P} \alpha_{m,i} \Delta t \mathcal{R}(q_{n-i}), \tag{35}$$

where $P$ is the number of previous points used and $M$ is the maximum order of the $\varphi$ function. Notice that if methods constructed using this ansatz are combined with adaptive Krylov-type algorithms (such as KIOPS) for calculating the linear combinations of terms like $\varphi_m(\Delta t \mathcal{J}_n) v_m$, then each time step will require only one call to KIOPS (or to any other adaptive-Krylov method). Since approximating the linear combinations involving $\varphi$ matrix functions-vector products represents the main computational cost of an EPI method, all such schemes will have a comparable cost. Therefore, significant increases in accuracy with only relatively small increases in computational cost are expected from such higher-order methods. This computational advantage of the newly introduced EPI schemes will hold for problems where the cost of the Krylov-projection evaluation of the exponential matrix functions and vectors dominates the computational cost per time step compared to the evaluation of the right-hand-side function or the actual product of the Jacobian and a vector. Should the latter computations become prohibitively expensive and dominant, one might still want to use a lower-order method.

To construct specific EPI schemes, the coefficients $\alpha_{m,i}$ in Eq. (35) must be determined. This is done with the Butcher trees order condition theory, which greatly simplifies their derivation. Details on Butcher trees to construct exponential methods are found in [58,57]. This machinery is used to derive systems of equations for the coefficients $\alpha_{m,i}$. In the case of multistep-type EPI methods, such systems are linear but typically under-determined. Therefore, as for standard Runge-Kutta and multistep methods, additional constraints may be added and the extended linear systems may be solved to obtain either a single method or a family of schemes. Such constraints may serve to introduce certain desired properties to the resulting scheme, for instance: 1) the largest possible number of zero coefficients, 2) coefficients with the smallest possible magnitude in the scheme itself or in the next order error term. The coefficients of a method may be represented as a matrix $A$, where $A_{i,j} = \alpha_{i,j}$. Therefore, a method using $P$ previously computed nodes $q_{n-i}$ ($i = 1, ..., P$) and $M$ functions $\varphi_j$ ($j = 1, ..., M$) will have $P$ columns and $M$ rows.

Here, it is chosen to supplement the order conditions by constraints that minimize the magnitude of the coefficients in the first term of the error. In other words, if the method is of order $p$, then the coefficients of the error term of order $p+1$ will be minimized. Solving the order conditions with these constraints, methods of orders 4, 5 and 6 are obtained with the following coefficients:

$$A_4 = \begin{pmatrix} 0 & 0 \\ -\frac{3}{10} & \frac{3}{40} \\ \frac{32}{5} & -\frac{11}{10} \end{pmatrix}, \tag{36}$$

$$A_5 = \begin{pmatrix} 0 & 0 & 0 \\ -\frac{4}{5} & \frac{2}{5} & -\frac{4}{45} \\ 12 & -\frac{9}{2} & \frac{8}{9} \\ 3 & 0 & -\frac{1}{3} \end{pmatrix}, \tag{37}$$

$$A_6 = \begin{pmatrix} 0 & 0 & 0 & 0 \\ -\frac{49}{60} & \frac{351}{560} & -\frac{359}{1260} & \frac{367}{6720} \\ \frac{92}{7} & -\frac{99}{14} & \frac{176}{63} & -\frac{1}{2} \\ \frac{485}{21} & -\frac{151}{14} & \frac{23}{9} & -\frac{31}{168} \end{pmatrix}. \tag{38}$$

Note that the EPI multistep methods require computation of the starting nodes $q_n$ to begin the time stepping process. These values may be calculated in many different ways. Here, the EPI2 method (see Eq. (25)) with a smaller time step to compute the necessary starting values is used. Once the initial values are obtained and the remainder function $R(q_n)$ is evaluated and stored, the time stepping algorithm proceeds as follows. At each new time iteration, the solution at the next time step $q_{n+1}$ is computed using Eq. (34) as well as the following steps:

1. compute vectors $v_m = \sum_{i=1}^{P} \alpha_{m,i} \mathcal{R}(q_{n-i})$ using previously stored values $q_{n-i}$ and $\mathcal{F}(q_{n-i})$. Note that while vectors $\mathcal{F}(q_{n-i})$ may be reused from previous iterations, vectors $\mathcal{R}(q_{n-i})$ cannot be reused because the Jacobian matrix changes from one iteration to the next;
2. use the KIOPS algorithm to evaluate

$$w_n = \varphi_1(\Delta t \mathcal{J}_n) \mathcal{F}(q_n) + \sum_{m=1}^{M} \varphi_m(\Delta t \mathcal{J}_n) v_m;$$

3. advance the solution over the next time step $q_{n+1} = q_n + w_n \Delta t$.

Exponential integrators, much like implicit methods, possess very good stability properties since they represent an exact solution $q_{n+1} = e^{\Delta t \mathcal{J}} q_n$ to the linear problem $y' = \mathcal{J} y$ (i.e., $\mathcal{R}(q) = 0$ and $\mathcal{J}_n = \mathcal{J}$ in Eq. (22)) with $\|e^{\Delta t \mathcal{J}_n}\| < 1$ if the real parts of the eigenvalues of $\mathcal{J}_n$ are all negative. In addition, exponential integration can lead to computational savings per time step compared to implicit schemes. This is the result of the Krylov methods performing more efficiently in evaluating exponential-like functions $\varphi$ in expression (29) compared to computing the rational functions of a Jacobian, which have to be approximated within an implicit method. Exponential methods are also generally more accurate than explicit and implicit methods of the same order. In the context of geophysical fluid dynamics, this means that relevant wave dispersion relations are expected to be simulated more accurately.

Note that all exponential integrators are trivially A-stable since an exponential method solves the linear part of the system exactly ($q_n + \varphi_1(\Delta t \mathcal{J}_n) \Delta t \mathcal{F}(q_n)$) thus no unstable computational modes are generated. Surely when the exponential is approximated rather than computed exactly the method is not necessarily A-stable. However, when adaptive Krylov-based algorithm is used to compute the products of exponential matrix functions and vectors these computations can be performed in a way that does not pose challenges for stability. It has been shown in the literature [59,60] that combining, exponential Runge-Kutta integration with Krylov algorithm is similar to employing an explicit Runge-Kutta method with coefficients that depend on the eigenvalues (particularly those from the boundary of the spectral domain) of the Jacobian and therefore the region of stability of such method also depends on these eigenvalues. Obviously, such explicit Runge-Kutta method will be different for each time step, as the Jacobian and the Krylov vectors are changing. Thus the stability region is not static and is accounting for the eigenvalues of the Jacobian. In addition, the adaptive Krylov allows to set tolerance on the approximation of the exponential matrix functions-vector products and it is easy to ensure that these computations are done with enough precision to not impact the overall time stepping process. It is also worth noting that a similar issue arises when implicit-Krylov methods are used and solving with an approximate Jacobian and within a certain tolerance is a common practice with implicit methods. The implications of this strategy in the context of Krylov subspace methods for approximating the matrix exponential have been studied recently in [53]. Practice shows that employing exponential-Krylov time stepping approach is effective and allows integration of many practical problems without issues with stability [61,62].

In the next section, numerical examples will be used to demonstrate that the newly constructed EPI methods allow significant increases in accuracy of the time integration without incurring significant computational cost.

## 5. Numerical experiments

The numerical properties of the algorithms described in previous sections are studied using test cases found in the literature. In §5.1, the accuracy and convergence of the EPI time integrators on three PDE problems are evaluated. In §5.2, numerical properties of the spatial discretization are evaluated in isolation using benchmarks found in [63]. In §5.3, the accuracy of the various time integrators and DFR accuracy orders is studied from simulating the barotropic test case found in [64]. Finally, results from three test cases based on the benchmarks found in [65] are presented. These experiments allow, among other things, to assess the conservation of various quantities.

Unless otherwise indicated, all tests with the shallow-water equations use a time step of 1 hour and the tolerance of the KIOPS solver is set to $10^{-10}$. The $\epsilon$ parameter of the complex-step approximation is set to $\approx 1.5 \times 10^{-8}$. The computational grid is made up of $6 \times N_e \times N_e$ elements. Each element consists of $N_s \times N_s$ solution points. The value of $N_s$ and $N_e$ are varied such that the total number of degrees of freedom is kept constant. For the considered shallow-water test cases, a global grid with 86400 degrees of freedom and a timestep of 1 hour corresponds to Courant numbers ranging from 10 to 20. The grid rotation parameters are set to $\lambda_0 = 0$, $\phi_0 = \pi/4$ and $\alpha_0 = 0$ (see Appendix B) since this configuration poses a greater challenge for most test cases than, say, the unrotated configuration.

Beside the intrinsic numerical diffusion of the Riemann solver, no filter, artificial diffusion or limiter is used in this study. These aspects have been extensively studied in the literature (see e.g. [66,6]) and may be assessed in future works.

Errors associated with simulated cases for which known analytical solutions exist are evaluated from the height field using the following global norms:

$$L_1(t) = \frac{I[|H - H_T|; t]}{I[H_T; t]}, \tag{39}$$

$$L_2(t) = \sqrt{\frac{I[(H - H_T)^2; t]}{I[H_T^2; t]}}, \tag{40}$$

$$L_\infty(t) = \frac{\max \left| H(x^1, x^2, t) - H_T(x^1, x^2, t) \right|}{\max \left| H_T(x^1, x^2, t) \right|}, \tag{41}$$

where $H_T$ is the height field of the analytical solution and $I$ is an approximation to a global integral calculated with the Gauß-Legendre quadrature rules as follows:

$$I[\Psi; t] \approx \sum_{p=0}^{5} \sum_{k=1}^{N_e} \sum_{l=1}^{N_e} \sum_{m=1}^{N_s} \sum_{n=1}^{N_s} \Psi\left((x^1)_{m,k,l,p}, (x^2)_{n,k,l,p}, t\right) \times$$
$$\sqrt{g}\left((x^1)_{m,k,l,p}, (x^2)_{n,k,l,p}, t\right) w_m w_n, \tag{42}$$

where $p$ is the panel index, the indices $k, l, m, n$ refer to the position of the points on the Gauß-Legendre grid and the $w$'s are quadrature weights that lead to an exact integration of polynomials of degree $2N_s - 1$ or less.

### 5.1. Convergence of the time integrators

The performance of the new EPI methods introduced above is evaluated with a standard set of tests for large stiff systems of equations [67]. Their convergence is assessed from three common benchmarks:

- Advection-diffusion-reaction PDE:

$$\frac{\partial u}{\partial t} + \alpha \left( \frac{\partial u}{\partial x} + \frac{\partial u}{\partial y} \right) = \varepsilon \left( \frac{\partial^2 u}{\partial x^2} + \frac{\partial^2 u}{\partial y^2} \right) + R(u),$$

where $R(u) = \gamma u(u - 1/2)(1 - u)$, $\varepsilon = 1/100, \alpha = -10, \gamma = 100$, on the domain $x, y \in [0, 1], t \in [0, 0.1]$ and with homogeneous Neumann boundary conditions. The initial condition is $u_0 = 256(xy(1-x)(1-y))^2 + 0.3$. The equation is discretized using a second-order finite difference method with $n = 1600$ grid points.

- Burger's equation:

$$\frac{\partial u}{\partial t} + \frac{1}{2}\frac{\partial u^2}{\partial x} = \varepsilon \frac{\partial^2 u}{\partial x^2},$$

with $\varepsilon = 10^{-3}$. This equation is discretized using a second-order finite difference method with $n = 1024$ grid points on the domain $x \in [0, 1], t \in [0, 1]$. The initial condition $u_0 = \exp\left(-(x-\mu)^2/(2\sigma^2)\right)$ with $\mu = 0.3, \sigma = 0.05$ and homogeneous Dirichlet boundary conditions are used.

- Semilinear parabolic PDE [67]:

$$\frac{\partial u}{\partial t} - \frac{\partial^2 u}{\partial x^2} = \int_0^1 u(x, s)ds + \phi(x, t),$$

where $\phi(x, t)$ is chosen such that the exact solution is $u(x, t) = x(1-x)e^t$. This equation is discretized in space using a second-order finite difference method for the derivative and the trapezoidal rule for the integral with $n = 400$ grid points on the domain $x \in [0, 1], t \in [0, 1]$. The initial condition $u_0 = x(1-x)$ and homogeneous Dirichlet boundary conditions are used.

The tolerance for KIOPS is set to $10^{-14}$ and the error is defined as the discrete 2-norm of the difference between the approximation to the solution and the reference solution. The exact Jacobian is used in all of the calculations. The reference solutions for the first two cases are computed using MATLAB's `ode15s` integrator with absolute and relative tolerances set to $10^{-14}$, whereas an exact solution exists for the third case.

Fig. 1 shows convergence plots obtained from the three test cases above. All methods converge to the reference solutions with the theoretically expected order of accuracy. Note that the semilinear parabolic PDE problem was constructed in [67] to demonstrate that certain time integration methods can suffer from the order reduction if they were not derived using the stiff order conditions. It has not been proven, however, that the stiff order conditions in [67,68] are necessary to avoid order reduction. While EPI schemes proposed here are not derived using these stiff order conditions from [67,68], no order reduction is observed for these schemes in the test problems.

Fig. 2 shows work-precision diagrams (CPU time vs. error) for all tested exponential scheme and benchmark cases. The points are in decreasing order of time step size, from left to right, with the following value of $h$: Advection-diffusion-reaction test case: $h = \{7.0 \times 10^{-4}, 2.2 \times 10^{-4}, 9.8 \times 10^{-5}\}$, Burger's equation test case: $h = \{3.1 \times 10^{-3}, 1.8 \times 10^{-3}, 1.0 \times 10^{-3}, 6.0 \times 10^{-4}, 3.0 \times 10^{-4}\}$ and semilinear parabolic test case: $h = \{1.1 \times 10^{-1}, 7.0 \times 10^{-2}, 5.0 \times 10^{-2}, 3.3 \times 10^{-2}, 2.1 \times 10^{-2}, 1.5 \times 10^{-2}, 1.0 \times 10^{-3}\}$.

The work-precision diagrams clearly show that increasing the order of the EPI scheme does not lead to a significant increase in the computational time required to approximate the solutions with specified accuracy. This is because all EPI methods are implemented using a single computation of a linear combination of $\varphi$ functions per time step and the computational cost of these linear combinations is comparable for all of the time integrators.

### 5.2. Convergence of the DFR method

The diffusion-free, zonally balanced, time-dependent flow proposed in [63] is a difficult test for the cubed-sphere when the grid is rotated at a 45° angle. The highest values of the height field cross eight panel edges, four corners of the cube and follow two panel edges located along the equator. This test case has an analytical solution and is therefore frequently used to assess the convergence of numerical models. All simulations presented in this subsection are integrated in time with the EPI6 scheme. Error norms are computed at day 5.

To illustrate the numerical convergence of the DFR method applied to the shallow-water equations, the experiments are organized in two ways. First, a $p$-convergence test is performed. The results are shown in Fig. 3a with $N_e = 10$ and the number of solution points is increased from $N_s = 3$ to $N_s = 8$. Then, an $h$-convergence test is conducted by varying $N_e$ from 10 to 20, while keeping a fixed number of solution points $N_s = 5$. The results are shown in Fig. 3b. Table 1 shows the order of accuracy calculated with respect to the errors with $N_e = 10$ and $N_e = 20$. It is observed that the DFR scheme reaches the expected formal order of accuracy, except at order 4 where there is a small order reduction.

### 5.3. Barotropic instability

A further description of convergence properties is obtained from a barotropic instability case using various numbers of elements and solution points. For this test proposed in [64], the initial condition is a zonal flow imitating a tropospheric jet at mid-latitudes in the northern hemisphere. A small perturbation of the height field is introduced to induce the development of barotropic instability. This test case describes the evolution of a barotropic wave, where a continuous transfer

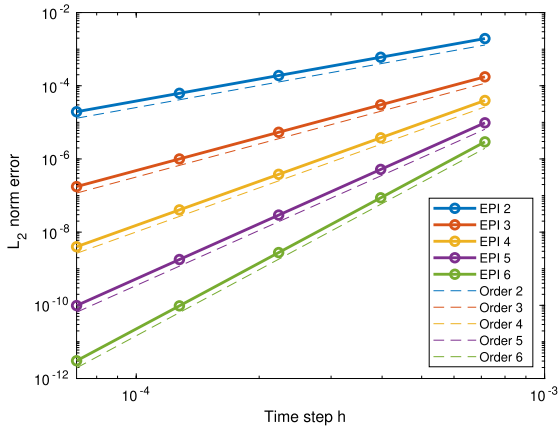
(a) Advection-diffusion-reaction

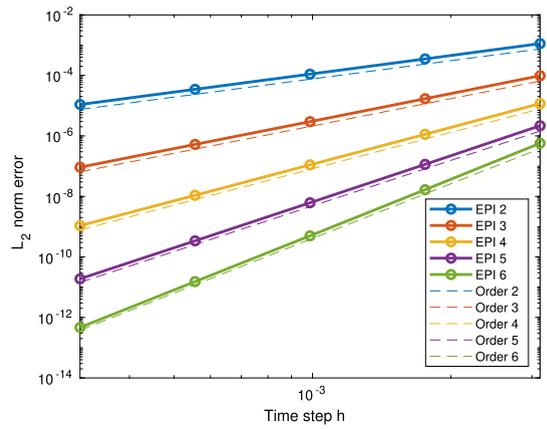
(b) Burger's equation

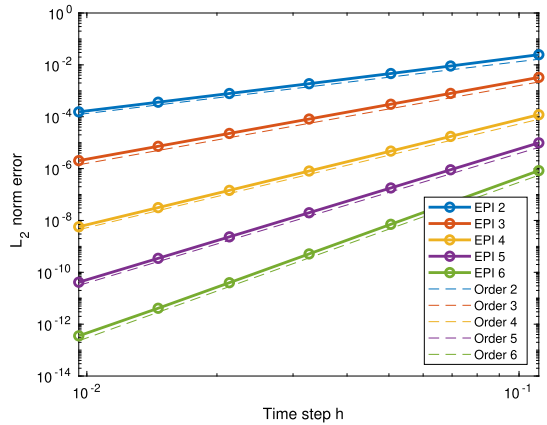
(c) Semilinear parabolic

**Fig. 1.** Convergence plots for multistep-type EPI methods of orders 2, 3, 4, 5 and 6.

**Table 1**
Computed order of accuracy for the diffusion-free, zonally balanced, time-dependent flow.

| $N_s$ | $L_1$ | $L_2$ | $L_\infty$ |
| --- | --- | --- | --- |
| 3 | 3.501 | 3.504 | 3.007 |
| 4 | 3.724 | 3.705 | 3.668 |
| 5 | 5.724 | 5.571 | 4.964 |
| 6 | 6.030 | 6.061 | 6.222 |

of energy occurs at different spatial scales at mid-latitudes. This test is challenging for non-monotonic numerical schemes, such as those presented in this study.

The relative vorticity

$$\zeta = \frac{\varepsilon^{0ij}}{\sqrt{g}} \frac{\partial}{\partial x^i} \left( g_{jk} u^k \right), \tag{43}$$

where $\varepsilon^{\alpha\mu\nu}$ is the Levi-Civita symbol, may be compared with the reference solution presented in [64]. Fig. 4 shows the relative vorticity for different values of $N_s$. The time integrator is chosen to match the order of convergence in space. The number of elements $N_e$ is changed accordingly in order to keep a constant number of degrees of freedom. Although the grid spacing is rather coarse, a convergence to the reference solution is observed as the order increases. The shapes of the solutions at 5th and 6th orders are comparable to the reference solution.

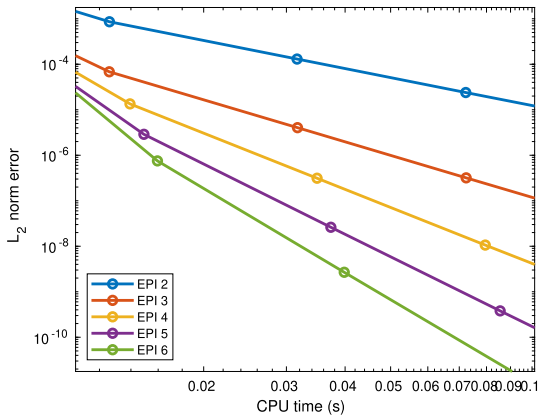
(a) Advection-diffusion-reaction

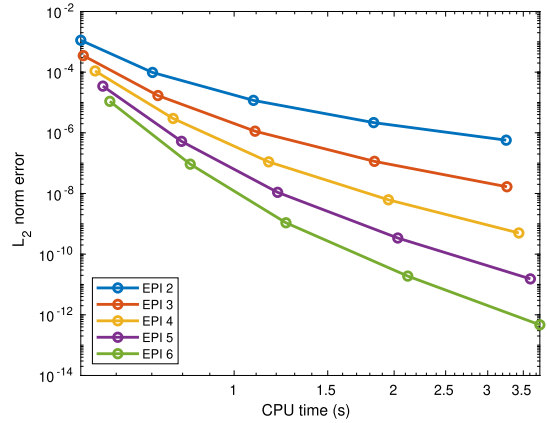
(b) Burger's equation

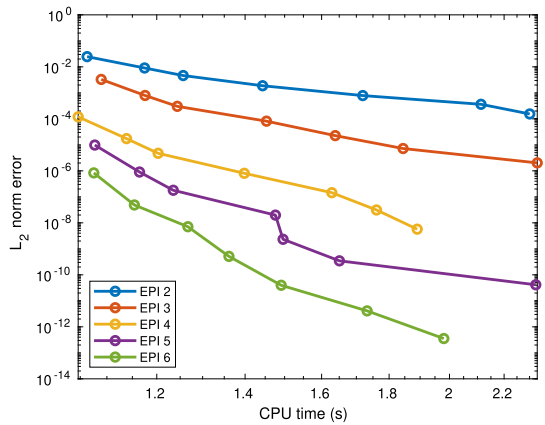
(c) Semilinear parabolic

**Fig. 2.** Work-precision diagrams for multistep-type EPI methods of orders 2, 3, 4, 5 and 6.

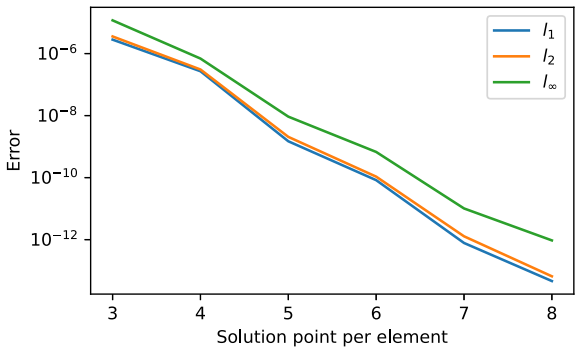
(a) $p$-convergence ($N_e = 10$)

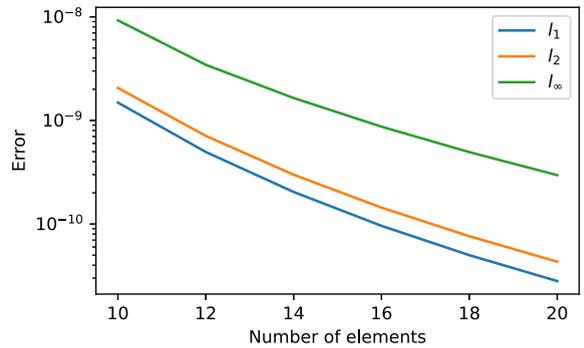
(b) $h$-convergence ($N_s = 5$)

**Fig. 3.** Error convergence at day 5 for the diffusion-free, zonally balanced, time-dependent flow.

### 5.4. Numerical conservation properties

Three standard tests suggested in [65] are considered to assess the conservation properties of the numerical schemes presented in the previous sections. To monitor the residual temporal evolution of global invariants, one defines a normalized conservation error as follows:

$$\hat{\Psi}(t) = \frac{I[\Psi; t] - I[\Psi; 0]}{I[\Psi; 0]}. \tag{44}$$

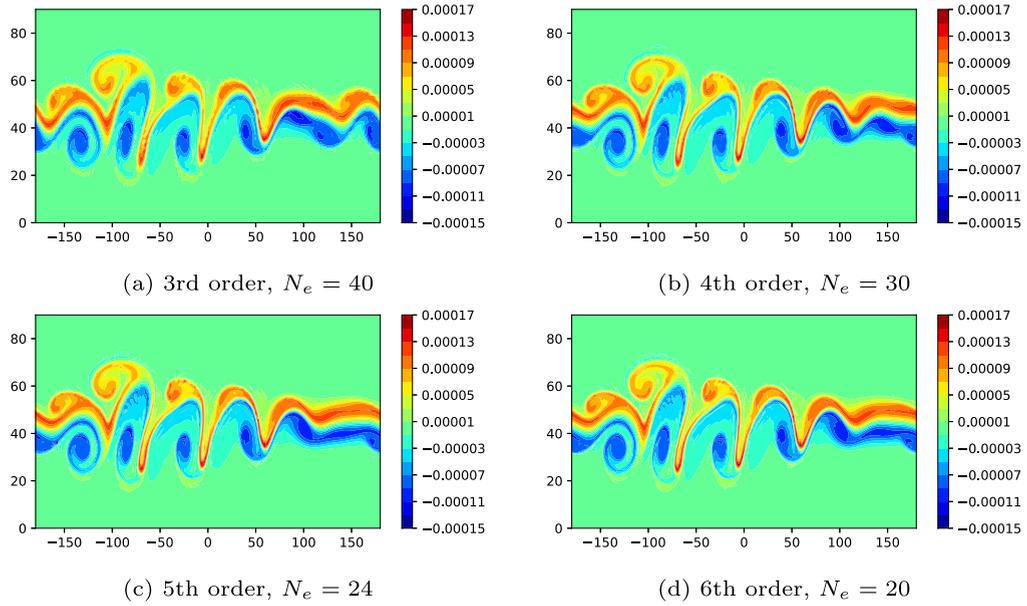

**Fig. 4.** Relative vorticity field associated with the barotropic instability test at day 6 on a global grid with 86400 degrees of freedom. Only the Northern Hemisphere is shown.

The function $\Psi$ is replaced by $H$ for mass conservation, by

$$\frac{1}{2}\left[Hg_{ij}u^i u^j + g_r\left((H+h_B)^2 - h_B^2\right)\right]$$

for total energy conservation, and by

$$\frac{(\zeta + f)^2}{2H}$$

for potential enstrophy conservation, where

$$f = \frac{\varepsilon^{0ij}}{\sqrt{g}} g_{jk}\Gamma^k_{i0} \tag{45}$$

$$= \frac{2\Omega}{\delta}(\sin\phi_p - X\cos\phi_p \sin\alpha_p + Y\cos\phi_p \cos\alpha_p) \tag{46}$$

is the Coriolis parameter.

In the following numerical experiments, four different configurations are considered. The value of $N_s$ varies from 3 to 6. The time integrator is chosen to match the order of convergence in space. The number of elements $N_e$ is changed in order to keep a constant number of degrees of freedom. The parameter $\alpha$ introduced in [65] is set to zero because the grid rotation mechanism in the present paper is embedded in the equations of motion.

The following approximation is proposed to compare the different mean resolutions of the cubed-sphere grids with Gauß-Legendre points to mean resolutions of other types of grids:

$$\Delta = \frac{90°}{N_e N_s}. \tag{47}$$

For the configurations considered in this subsection, $\Delta \approx 1°$. This resolution is comparable to the resolution of the models presented in [69].

#### 5.4.1. Steady-state geostrophically balanced flow

This test case is a steady-state solution to the shallow-water equations. The winds are a solid-body rotation, and the height is defined such that an exact geostrophic balance exists. The solution is expected to maintain this steady-state balance for at least 5 days.

Fig. 5 depicts the normalized error of the height field for the steady-state geostrophically balanced flow after 5 days. As the order increases, the error rapidly decreases. All configurations maintain a balance between advection, fluid height gradient, and Coriolis terms. Contrary to results found in other studies with the Yin-Yang overset grid [69,70] in which relatively large numerical errors occur where the zonal flow crosses the boundary between two panels, cubed-sphere panel boundaries and orientations do not impede convergence for higher-order configurations. The error associated with orders higher than 4 is smaller than the errors of the four models in [69], including the reference spectral model.

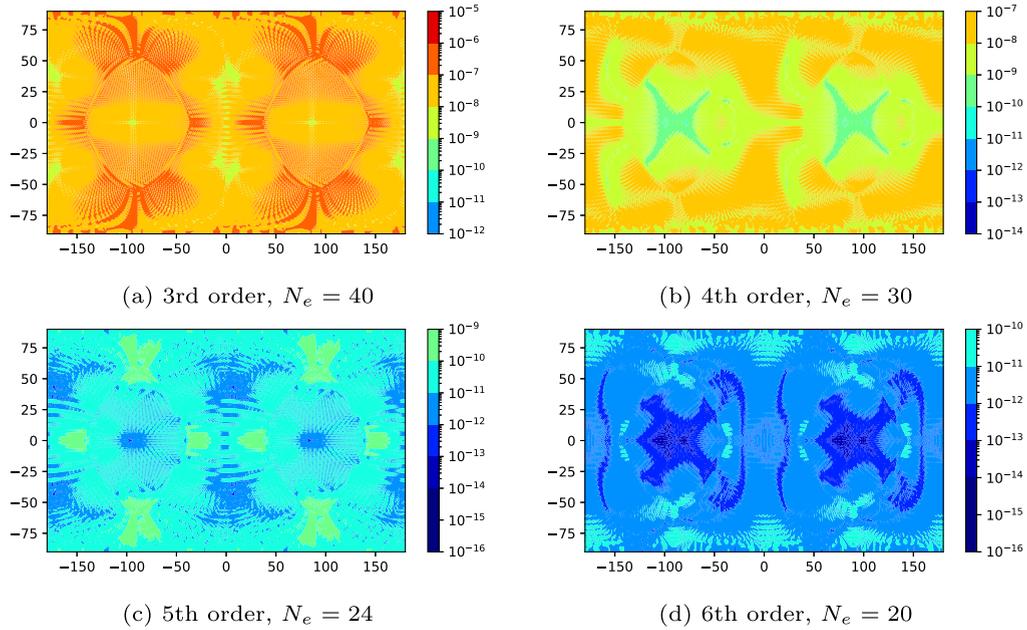

**Fig. 5.** Normalized error of the height field for the steady-state geostrophically balanced flow after 5 days using orders 3, 4, 5 and 6. A global grid with 86400 degrees of freedom is used.

Fig. 6 shows the time trace of the normalized error for the height field over a period of 30 days. The evolution of the normalized conservation error of the mass, total energy and potential enstrophy is presented in Fig. 7. Adjustments are observed during the first few time steps, likely due to imperfections in the initial conditions, but the conservation errors remain reasonably small for the rest of the simulations. Further experiments will be required to fully understand the convergence during the first few time steps. It is worth mentioning that nothing is done here to ensure that the discrete initial condition is in geostrophic balance. An example of a procedure to numerically enforce the geostrophic equilibrium is proposed in [71]. However, the numerical stability does not seem to be impacted by these initial imperfections since the error stabilizes quickly and remains low for the rest of the simulations.

*5.4.2. Zonal flow over an isolated mountain*

In this test case, a flow is perturbed by a topographically induced source term. The mountain shape is cone-like and not differentiable. In a model based on DFR, one must therefore approximate the mountain shape with piecewise polynomials, which may introduce spurious oscillations. In addition, the wind and height fields would be initially balanced only in the absence of a mountain. This results in a flow evolution that is difficult to predict numerically.

Orographic Rossby waves are induced from the beginning of the simulation. These waves then spread over almost the entire sphere, including the southern hemisphere. The height fields after 15 days are shown in Fig. 8. These results are similar to those reported by other authors (for instance [25,69]). There are no apparent non-physical features such as oscillations around the mountain, although there are visible differences in the shape of the solution of the different configurations. In order to better distinguish the differences between the simulations, the relative vorticity $\zeta$ is shown in Fig. 9. These results are comparable to those presented in [27]. Small oscillations are visible in the relative vorticity field, especially around the mountain. This could be caused by the approximation of the mountain slope by piecewise polynomials or could be a sign that aliasing errors are present in higher-order simulations. The use of filters and other strategies against aliasing-driven errors will be studied in future works. Chapter 5 of [29] provides useful informations on this topic.

The normalized errors of mass, total energy and potential enstrophy conservation over 30 days are shown in Fig. 10. The conservation of mass is accurate to machine precision and is comparable to results presented in [22]. The errors on the conservation of total energy and potential enstrophy are larger because the discretization of Eqs. (1) and (2) does not explicitly enforce their conservation.

It should also be mentioned that some of the errors may come from the fact that this study focuses on dynamical aspects relevant to atmospheric models. For instance, the formulation and space-time discretization presented here do not maintain a perfect equilibrium for some special cases such as lakes at rest (see for example [72]). Also, recall that nothing is done here to ensure that the discrete initial condition are in geostrophic balance. Thus, the restoring forces seek to compensate for these imbalances and other errors due to spatial discretization, leading to oscillations in the conservation error curves. This is particularly visible for energy and potential enstrophy, which are not conserved fields in the discrete system.

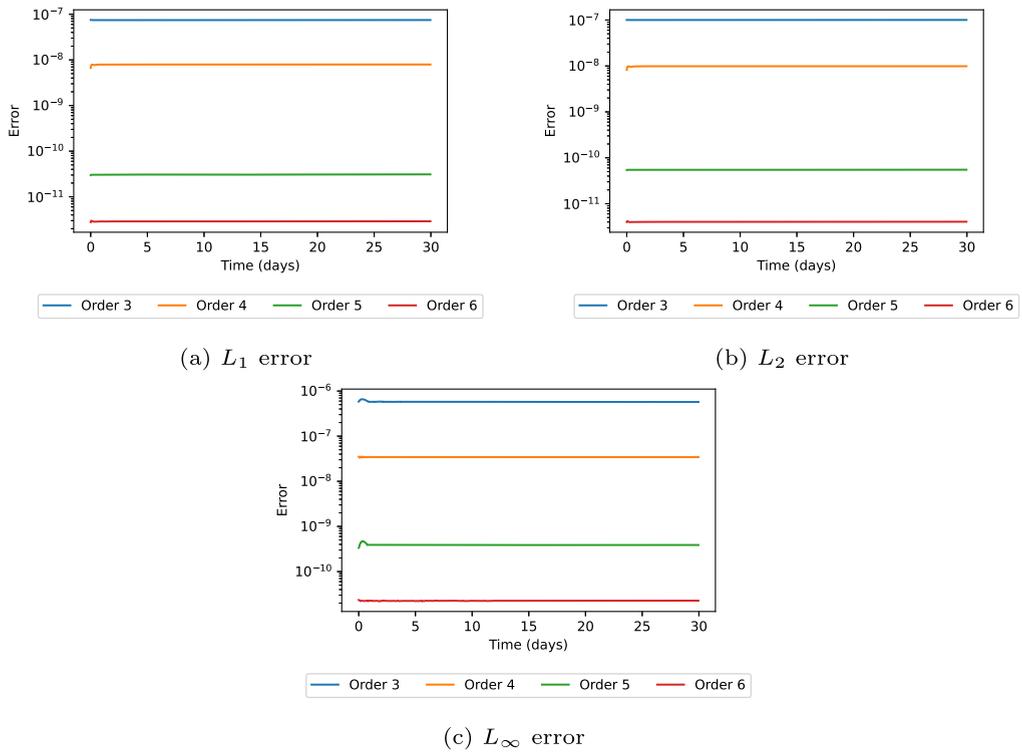

Fig. 6. Normalized error of the height field as a function of time for the steady-state geostrophically balanced flow.

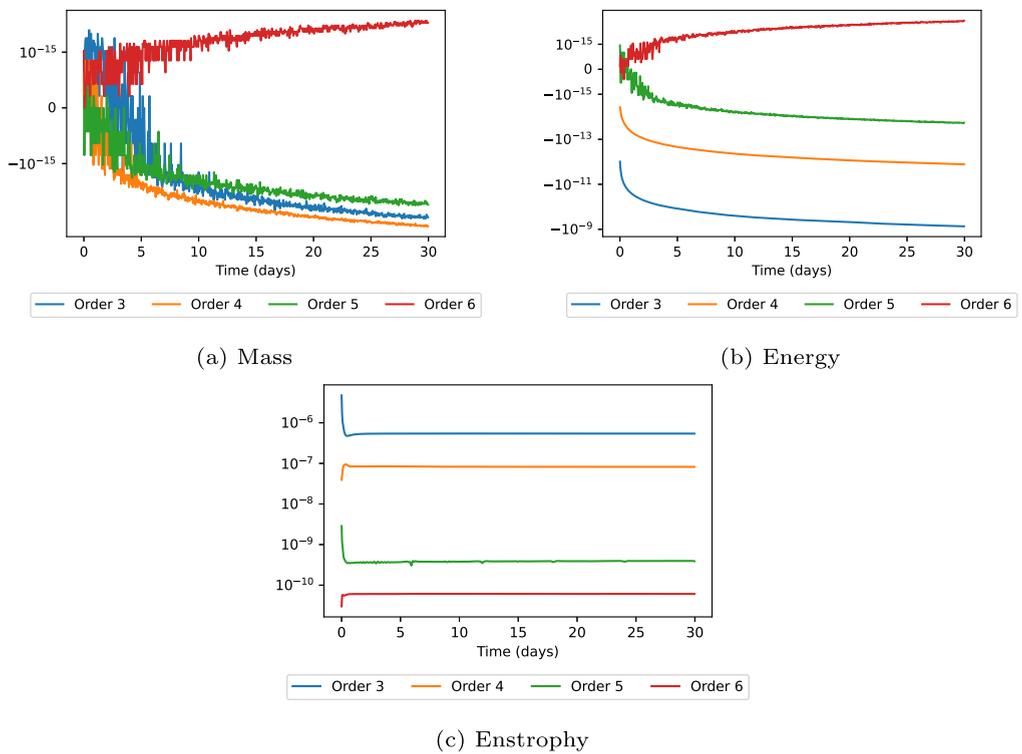

Fig. 7. Time traces of the normalized errors of conserved quantities for the steady-state geostrophically balanced flow.

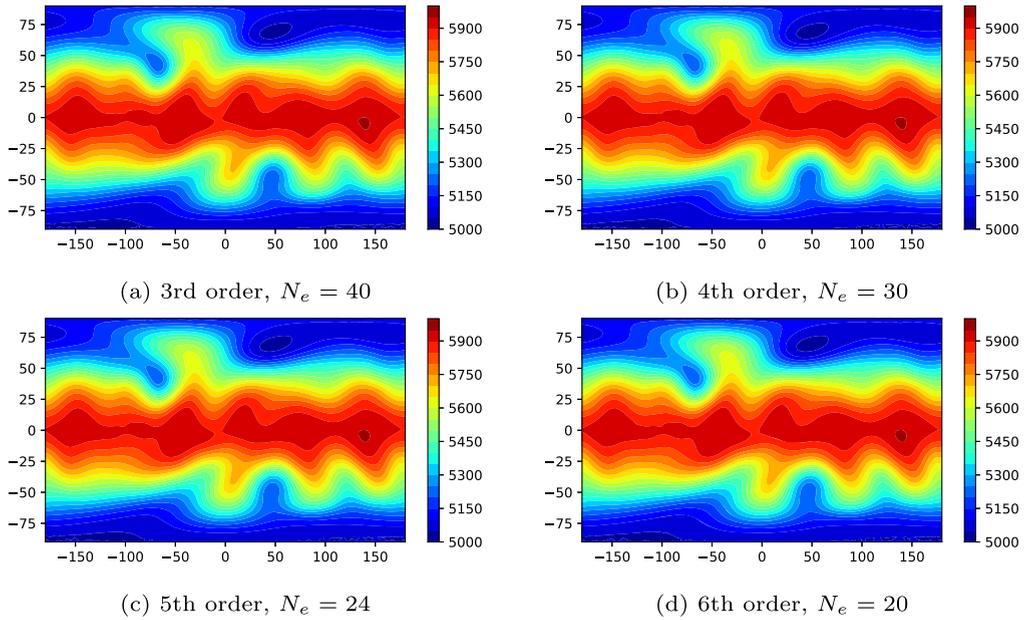

**Fig. 8.** Height field of the zonal flow over an isolated mountain at day 15 using orders 3, 4, 5 and 6. A global grid with 86400 degrees of freedom is used.

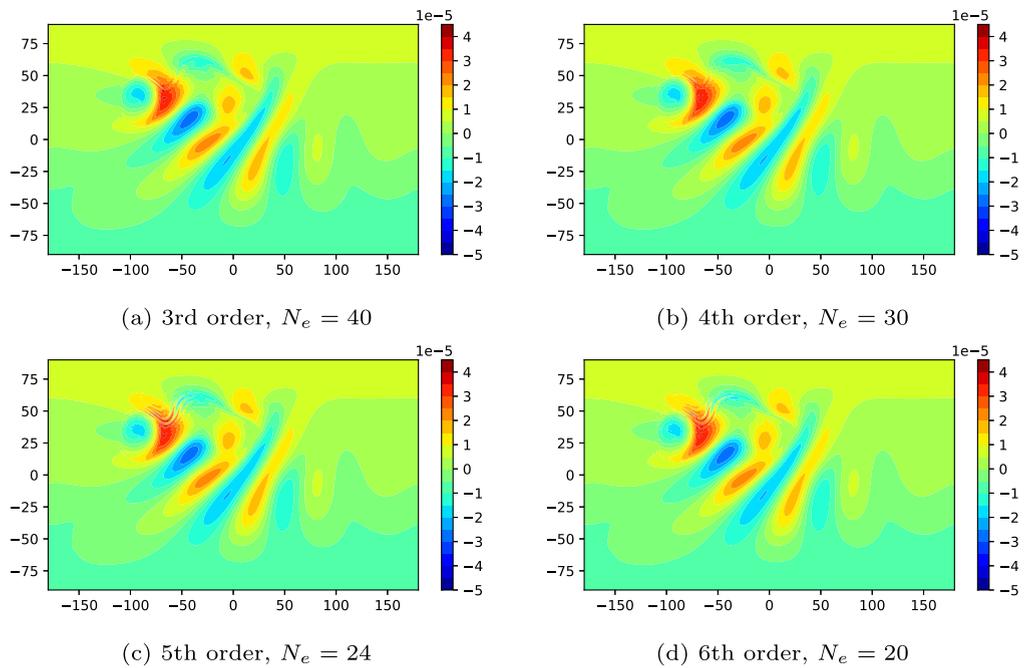

**Fig. 9.** Relative vorticity field of the zonal flow over an isolated mountain at day 7 using orders 3, 4, 5 and 6. A global grid with 86400 degrees of freedom is used.

#### 5.4.3. Rossby-Haurwitz wave

In this test case, the Rossby-Haurwitz wave number 4 is considered. This wave is an analytical solution to the non-linear barotropic vorticity equation. It is well known that this test case is susceptible to instabilities due to truncation in the initial conditions [73] and that the numerical solution will eventually lose its structure. However, the solution is expected to remain stable over the 14 days required by [65]. The wave number 4 is expected to propagate steadily and retain its structure with only slight wavering in its shape.

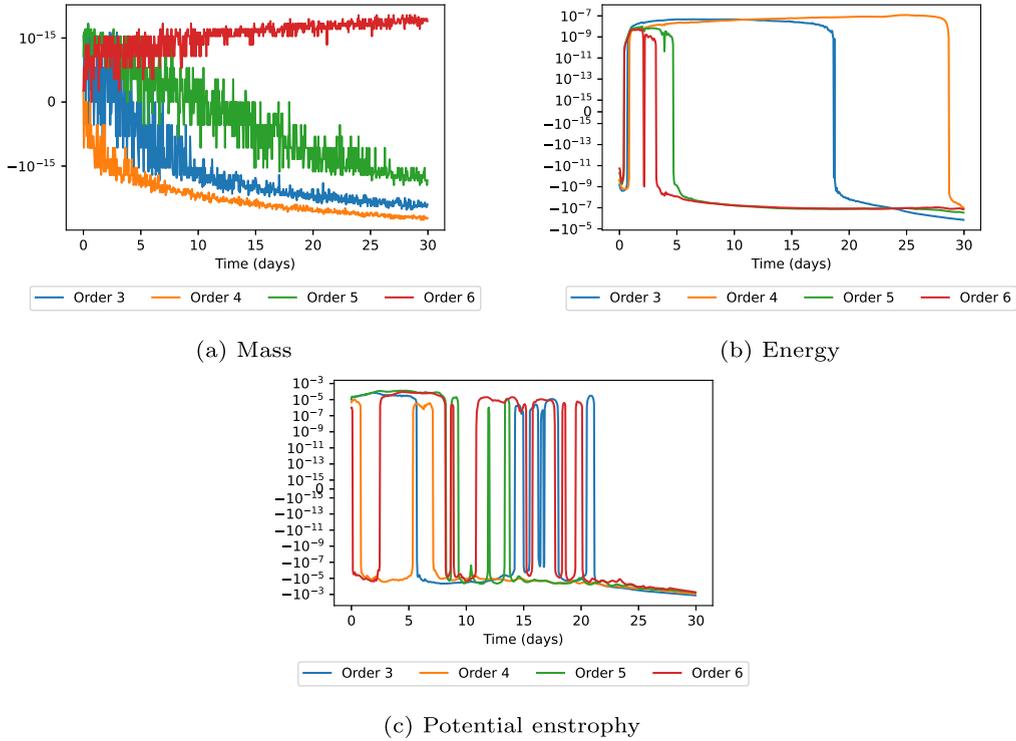

(a) Mass

(b) Energy

(c) Potential enstrophy

**Fig. 10.** Time traces of the normalized errors of conserved quantities for the zonal flow over an isolated mountain.

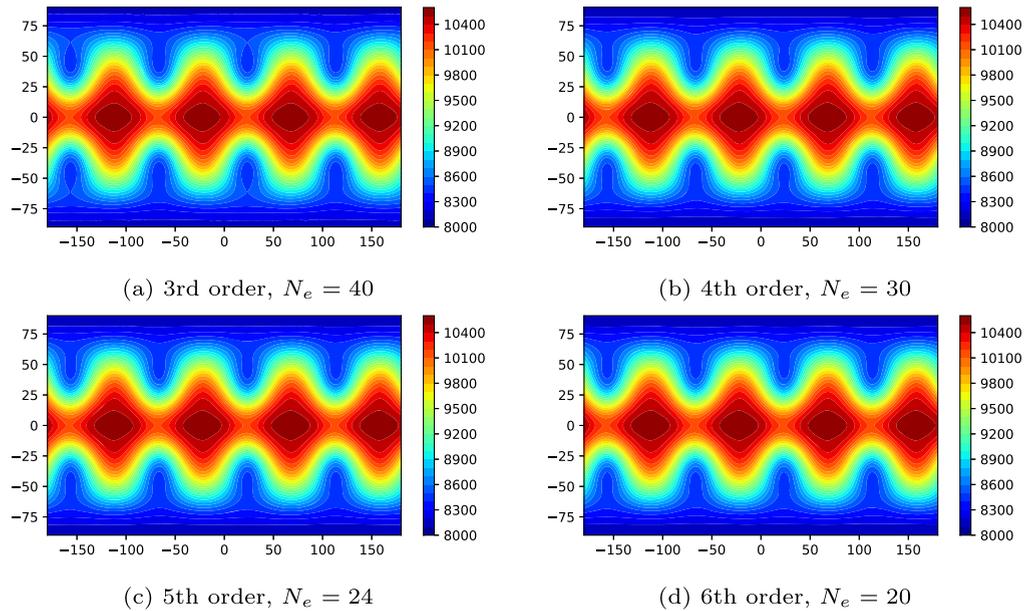

(a) 3rd order, $N_e = 40$

(b) 4th order, $N_e = 30$

(c) 5th order, $N_e = 24$

(d) 6th order, $N_e = 20$

**Fig. 11.** Height field of Rossby-Haurwitz wave at day 14 using orders 3, 4, 5 and 6. A global grid with 86400 degrees of freedom is used.

The simulated height fields after 14 days are shown in Fig. 11. The numerical solutions reproduce the reference solution reported in the literature [69] with good accuracy. In particular, the Rossby-Haurwitz wave remains stable throughout the duration of the simulation. The history plots of the conservation errors for mass, total energy and potential enstrophy are shown in Fig. 12. Mass is well conserved as in previous test cases. Conservation errors are much more important during the first 5 days of the simulation. Afterwards, these errors stabilize around a small value.

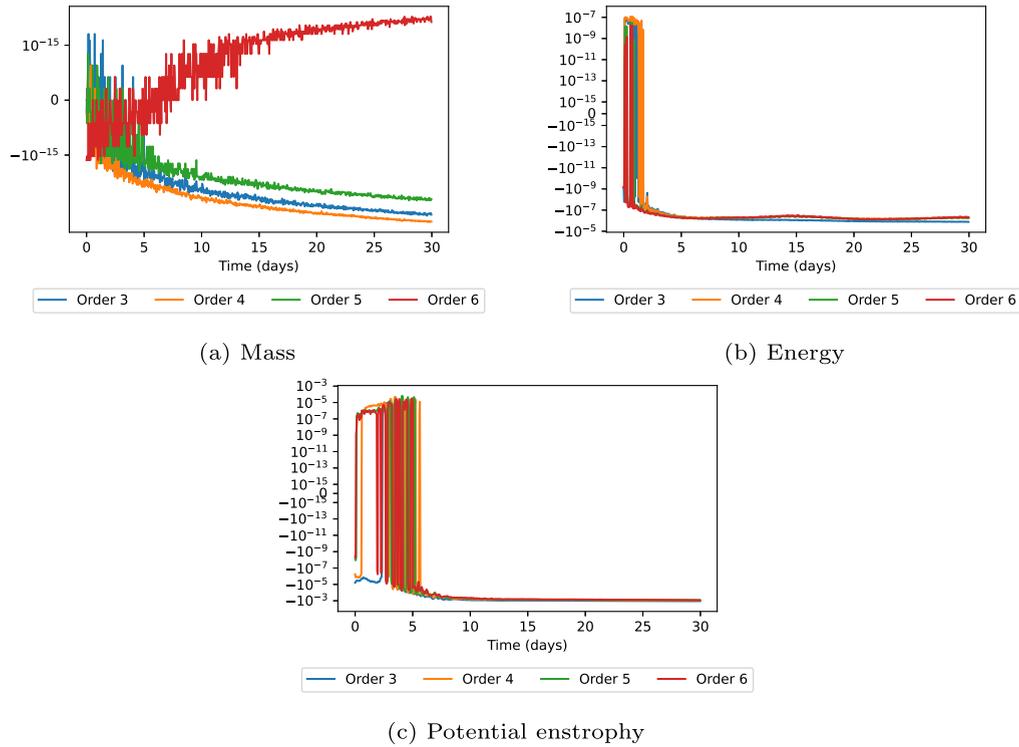

Fig. 12. Time traces of the normalized errors of conserved quantities for the Rossby-Haurwitz wave.

*5.5. Large time step*

As mentioned above, an important advantage of the exponential time integrators is that they allow large time steps regardless of the CFL condition. So far in this work, a time step of 1 hour has been used, which is larger than those typically allowed by an explicit Runge-Kutta method. Figs. 13 and 15 show the results of the zonal flow over an isolated mountain and the Rossby-Haurwitz wave with a large time step of 4 hours. For these configurations, the Courant number varies from 40 to 80. These numerical solutions are similar to those of Figs. 8 and 11. The differences between the results with a time step of 4 hours and 1 hour are shown in Figs. 14 and 16. On the one hand, for the zonal flow over an isolated mountain, the solution with a time step of 4 hours seems to converge as the order is increased. On the other hand, in the case of the Rossby-Haurwitz wave, significant differences are visible in the solutions of orders 5 and 6, while the solutions of orders 3 and 4 seem less sensitive to the increase in the time step size. This could be another indication that aliasing errors are present in high-order solutions. The success of these simulations can be largely attributed to the high accuracy with which the linear part of the equations is solved.

## 6. Summary and conclusions

The high-order methods discussed in this paper were applied to the shallow-water equations on the rotated cubed-sphere grid. Results from various test cases indicate that the proposed numerical algorithms produce realistic results. In particular, the convergence of the spatial and temporal patterns is well behaved as the order of accuracy is increased. Results are comparable to those obtained from more complicated methods based on the integral form of the equations of motion. Experiments show that mass is conserved at machine precision. A notable result is the ability of the schemes to perform accurate simulations with very large time steps for flows with complex structures.

The EPI integrators presented in this work offer not only high accuracy but new pathways to further improvements in making exponential integrators more computationally appealing. The authors plan to pursue this line of research in their future work. An important task ahead is the extension of this work to a comprehensive three-dimensional atmospheric model with the space-time tensorial approach to describe the Euler equations. Improvements to the performance of the exponential solver is underway to bring the algorithm closer to the performance obtained with semi-implicit semi-Lagrangian methods [74].

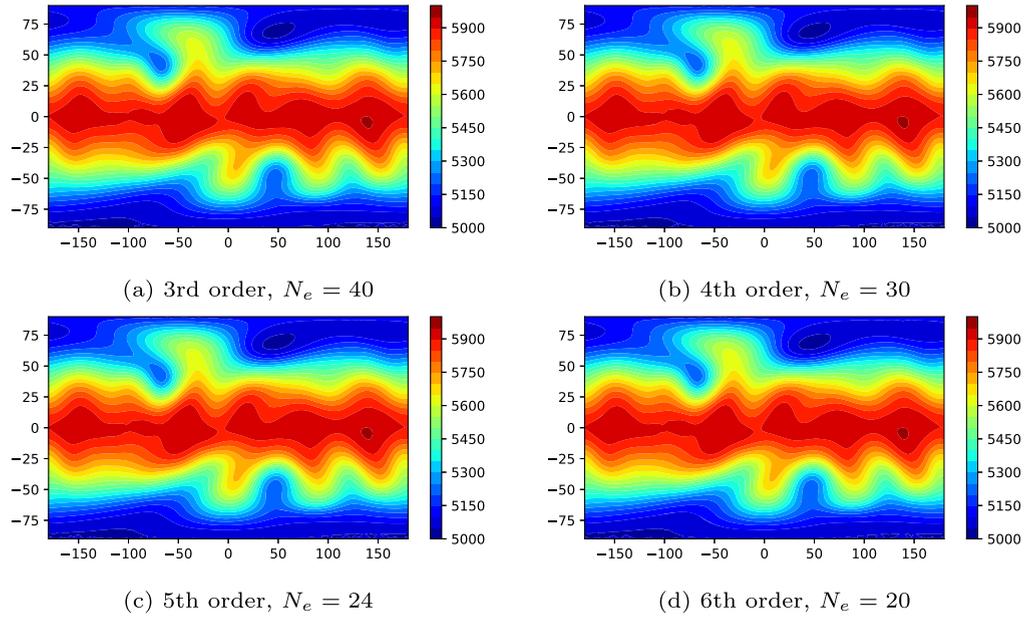

**Fig. 13.** Height field of the zonal flow over an isolated mountain at day 15 using orders 3, 4, 5 and 6 with a large timestep size of 4 hours. A global grid with 86400 degrees of freedom is used.

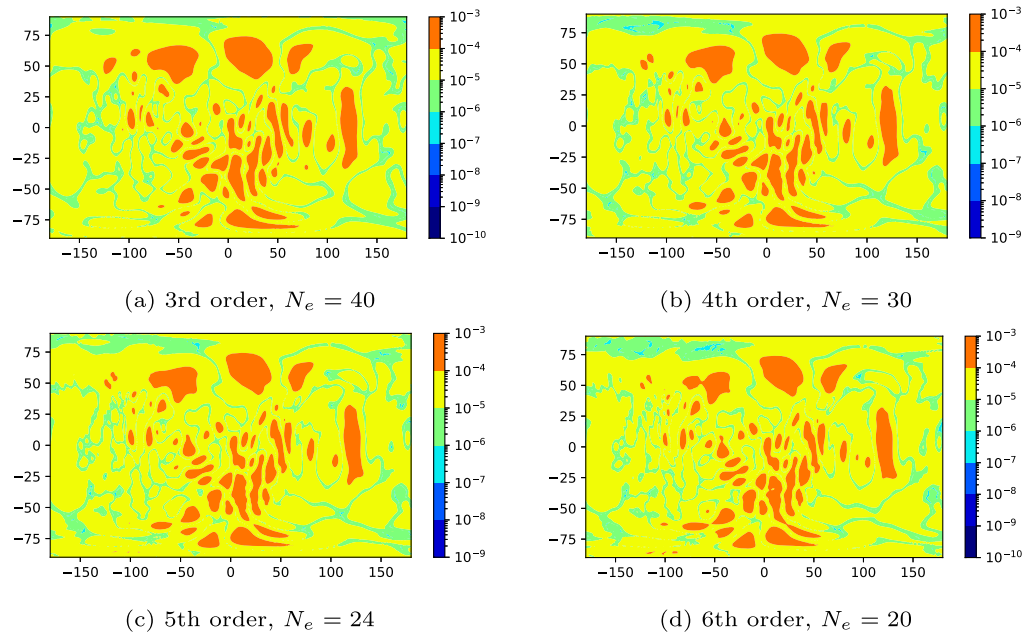

**Fig. 14.** Difference between timestep sizes of 4 hours and 1 hour for the zonal flow over an isolated mountain at day 15 using orders 3, 4, 5 and 6. A global grid with 86400 degrees of freedom is used.

## 7. Code availability

The code used to generate the results in this work is available under the GNU Lesser General Public License (LGPL) version 2.1 at https://doi.org/10.5281/zenodo.5014876 and https://gitlab.com/stephane.gaudreault/jcp2021_highorder_sw.

The EPIC package implements the exponential integrators and test cases presented in §5.1. The code and its license are available at http://faculty.ucmerced.edu/mtokman/#software.

A MATLAB implementation of KIOPS is available under the GNU LGPL version 2.1 at https://gitlab.com/stephane.gaudreault/kiops.

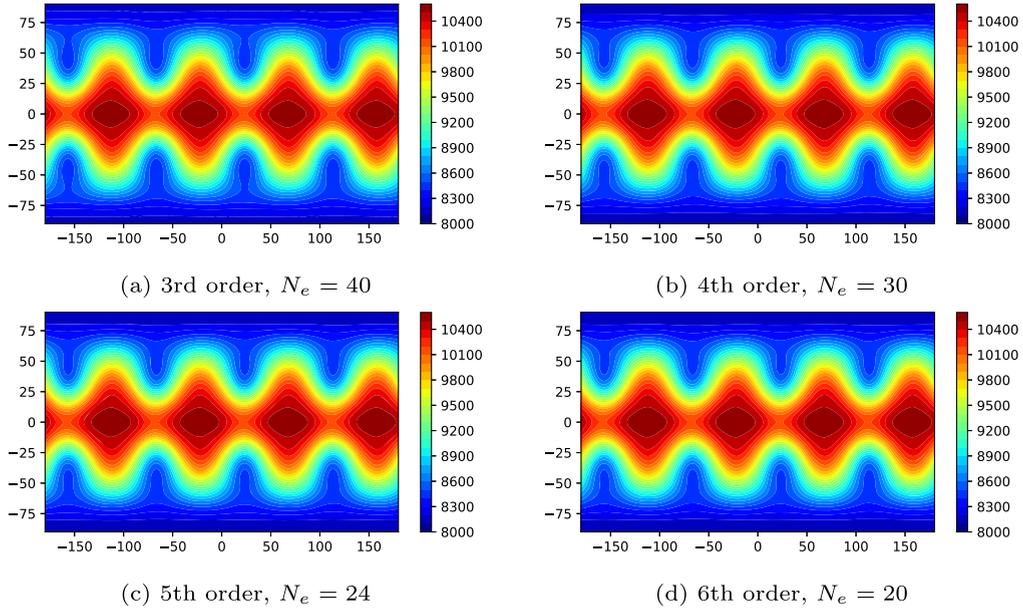

**Fig. 15.** Height field of the Rossby-Haurwitz wave at day 14 using orders 3, 4, 5 and 6 with a large timestep size of 4 hours. A global grid with 86400 degrees of freedom is used.

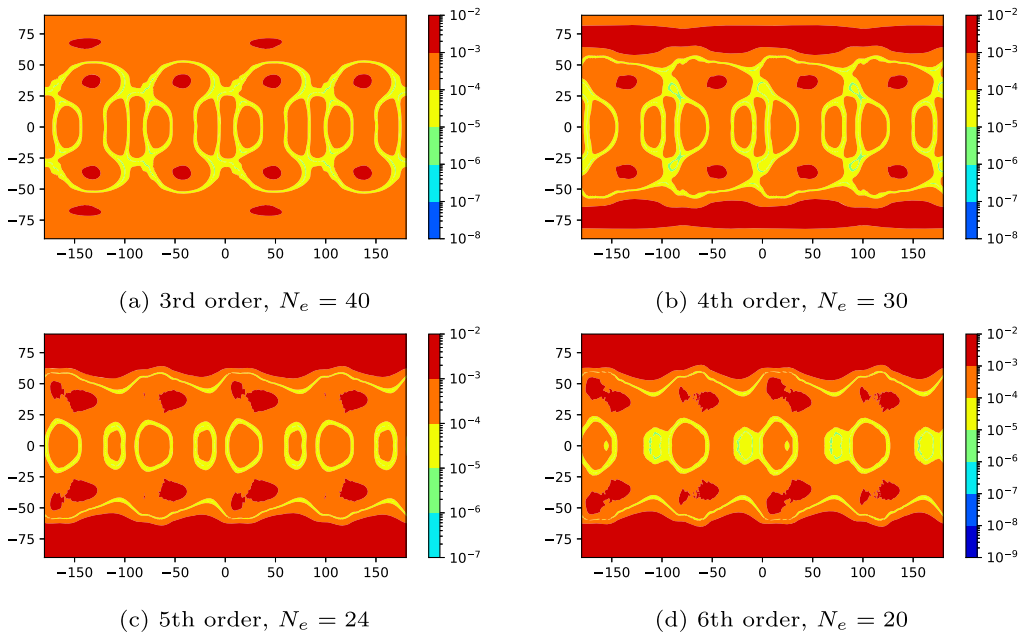

**Fig. 16.** Difference between timestep sizes of 4 hours and 1 hour for the Rossby-Haurwitz wave at day 14 using orders 3, 4, 5 and 6. A global grid with 86400 degrees of freedom is used.

**CRediT authorship contribution statement**

**Stéphane Gaudreault:** Conceptualization, Software, Validation, Writing, Project management. **Martin Charron:** Conceptualization, Software, Writing. **Valentin Dallerit:** Conceptualization, Software, Validation, Writing. **Mayya Tokman:** Conceptualization, Writing

**Declaration of competing interest**

The authors declare that they have no known competing financial interests or personal relationships that could have appeared to influence the work reported in this paper.

**Acknowledgements**

The authors sincerely thank Christopher Subich, Ayrton Zadra and Janusz Pudykiewicz for providing valuable comments on an earlier draft of this manuscript. They also thank two anonymous reviewers and Pedro Peixoto for providing comments on this manuscript. This work was supported in part by a grant no. 1115978 from the National Science Foundation, Computational Mathematics Program.

**Appendix A. The covariant shallow-water equations in 2+1 dimensions**

Here, the (2+1)-dimensional shallow-water equations will be derived from the (3+1)-dimensional covariant form of the Euler equations when the fluid is externally forced by a gravitational potential $\Phi$. Following [75],

$$\frac{\partial}{\partial t}\left(\sqrt{g}\rho\right) + \frac{\partial}{\partial x^j}\left(\sqrt{g}\rho u^j\right) = 0, \tag{48}$$

$$\frac{\partial u^i}{\partial t} + u^j \frac{\partial u^i}{\partial x^j} + \Gamma^i_{00} + 2\Gamma^i_{j0} u^j + \Gamma^i_{jk} u^j u^k = -h^{ij}\left(\frac{1}{\rho}\frac{\partial p}{\partial x^j} + \frac{\partial \Phi}{\partial x^j}\right) \tag{49}$$

describe the continuity and inviscid momentum equations ($i = 1, 2, 3$), respectively. Repeated Latin (spatial) indices are summed from 1 to 3, unless otherwise indicated. These equations hold under the classical assumptions of Newtonian space-time: time intervals are absolute, and space is also absolute. The symbol $\rho$ represents the fluid's density field, $p$ the pressure field, and $u^i \equiv dx^i/dt$ the velocity field components.

One may define an infinitesimal space interval $dl$ in an inertial frame as $dl^2 = h_{\mu\nu} dx^\mu dx^\nu$, where repeated Greek indices are summed from 0 to 3 (0 being a time index), and an infinitesimal time interval $dx^0 = u^0 dt$ as $(dx^0)^2 = \delta^0_\mu \delta^0_\nu dx^\mu dx^\nu$, where $\delta^\nu_\mu$ is the Kronecker tensor and $u^0$ an arbitrary non-zero constant with units of velocity. The time unit may be chosen such that $u^0 = 1$ without loss of generality and without impacting the form of the equations. Therefore, $u^0 = 1$ is assumed from now on. A (3+1)-dimensional distance $ds$ is written $ds^2 = dl^2 + (dx^0)^2 = g_{\mu\nu} dx^\mu dx^\nu = (h_{\mu\nu} + \delta^0_\mu \delta^0_\nu) dx^\mu dx^\nu$. Therefore, the contravariant tensor $h^{\mu\nu}$ is obtained from the (3+1)-dimensional contravariant metric tensor $g^{\mu\nu}$ (where $g^{\mu\nu} g_{\nu\alpha} = \delta^\mu_\alpha$):

$$h^{\mu\nu} = g^{\mu\alpha} g^{\beta\nu} h_{\alpha\beta} = g^{\mu\alpha} g^{\beta\nu} (g_{\alpha\beta} - \delta^0_\alpha \delta^0_\beta) = g^{\mu\nu} - g^{0\mu} g^{0\nu},$$

with $g^{00} = 1$ in Newtonian mechanics. This definition implies that $h^{\mu 0} = h^{0\mu}$ vanishes. Notice that in Newtonian space-time, there are in fact two tensors describing the geometry (for instance, $g_{\mu\nu}$ for space-time and $h^{\mu\nu}$ for space-only) because both spatial distances $dl$ and temporal intervals $dx^0$ are absolute and invariant.

The $\Gamma$'s in Eq. (49) are Christoffel symbols of the second kind: $\Gamma^i_{00}$ corresponds to the centripetal acceleration, $\Gamma^i_{j0}$ is associated with the local Coriolis acceleration, and $\Gamma^i_{jk}$ with the nonlinear metric terms. They are obtained from the usual definition in a metric space-time:

$$\Gamma^\mu_{\alpha\beta} = \frac{1}{2} g^{\mu\nu} \left(\frac{\partial g_{\nu\alpha}}{\partial x^\beta} + \frac{\partial g_{\nu\beta}}{\partial x^\alpha} - \frac{\partial g_{\alpha\beta}}{\partial x^\nu}\right) = \frac{1}{2} h^{\mu\nu} \left(\frac{\partial g_{\nu\alpha}}{\partial x^\beta} + \frac{\partial g_{\nu\beta}}{\partial x^\alpha} - \frac{\partial g_{\alpha\beta}}{\partial x^\nu}\right). \tag{50}$$

The last equality in Eq. (50) stems from the fact that $\Gamma^0_{\alpha\beta}$ vanishes in Newtonian mechanics. The term $\sqrt{g}$ is the square root of the determinant of the covariant metric tensor $g_{\mu\nu}$. In Newtonian mechanics, it also corresponds to the inverse of the square root of the 3-dimensional determinant of $h^{ij}$. Details on this (3+1)-dimensional formalism are provided in [75].

Consider spherical geometry under the thin-shell (or shallow-atmosphere) approximation. The $x^3$ axis is taken to represent the radial direction from the origin of the coordinate system and is chosen to indicate a geometric distance from the center of the Earth. The contravariant metric tensor components $g^{\mu\nu}$—and therefore $g_{\mu\nu}$ and $h^{ij}$—become independent of the coordinate $x^3$ because $x^3 = a$, where $a$ is the mean radius of the Earth in spherical geometry, as explained in [75,76]. The Christoffel symbols must not be directly approximated using $x^3 = a$ but rather recalculated from Eq. (50) with the approximated metric tensors (this preserves conservation laws). It turns out that this choice of coordinates together with the thin-shell approximation implies that the components $\Gamma^i_{\mu 3}$ of the Christoffel symbols vanish, as shown in [75] (their Eqs. (117)–(118)), and that $h^{13} = h^{23} = 0$ and $h^{33} = 1$ (their Eq. (60)).

Applying the spherical geopotential approximation in a geopotential coordinate, it may be shown that

$$h^{ij} \frac{\partial \Phi}{\partial x^j} + \Gamma^i_{00} = \begin{cases} 0 & \text{for } i = 1, 2; \\ g_r & \text{for } i = 3, \end{cases} \tag{51}$$

where $g_r \approx 9.80616$ m s$^{-2}$ is the constant effective gravitational acceleration at the surface of the Earth.

The "horizontal" (i.e. any spatial direction perpendicular to $x^3$) momentum equations become

$$\frac{\partial u^i}{\partial t} + u^j \frac{\partial u^i}{\partial x^j} + 2\Gamma^i_{j0} u^j + \Gamma^i_{jk} u^j u^k = -\frac{1}{\rho} h^{ij} \frac{\partial p}{\partial x^j} \tag{52}$$

for $i = 1, 2$, while the radial momentum equation leads to the hydrostatic balance:

$$\frac{\partial p}{\partial x^3} = -\rho g_r. \tag{53}$$

Recall that in Eqs. (48) and (52), all the metric terms (the $\Gamma$'s, $h^{ij}$ and $\sqrt{g}$) are independent of $x^3$ because the thin-shell approximation is being used.

To obtain the set of shallow-water equations, one makes the additional assumption that the fluid's density $\rho$ is uniform and constant. Integrating Eq. (48) over $x^3$ from $x^3 = h_B(t, x^1, x^2)$, where $h_B$ is prescribed, to $x^3 = h(t, x^1, x^2)$ while assuming that "horizontal" motion does not vary significantly over the depth $h - h_B$, one obtains

$$(h - h_B) \frac{\partial (\sqrt{g})}{\partial t} + \sqrt{g} \left( u^3(h) - u^3(h_B) \right) + (h - h_B) \frac{\partial}{\partial x^j} \left( \sqrt{g} u^j \right) = 0,$$

where in this case $j$ is summed from 1 to 2. Writing $u^3(h) = dh/dt$ and $u^3(h_B) = dh_B/dt$, and expanding the operator $d/dt = \partial/\partial t + u^j \partial/\partial x^j$ (where $j$ is summed from 1 to 2 because $h$ is independent of $x^3$), one gets the shallow-water continuity equation in arbitrary "horizontal" coordinates:

$$\frac{\partial}{\partial t} \left( \sqrt{g} H \right) + \frac{\partial}{\partial x^j} \left( \sqrt{g} H u^j \right) = 0, \tag{54}$$

where $H \equiv h - h_B$ is the thickness of the incompressible fluid and $j$ is summed from 1 to 2.

The 2-dimensional shallow-water momentum equations follow from assuming that the pressure at $x^3 = h(t, x^1, x^2)$ is a constant $p_0$. Therefore from Eq. (53), $p = \rho g_r (h - x^3) + p_0$, and the "horizontal" pressure gradient in Eq. (52) becomes $\partial p/\partial x^j = \rho g_r \partial h/\partial x^j$ ($j = 1, 2$).

A quasi-flux form may be obtained by adding Eq. (52) multiplied by $\sqrt{g} H$ and Eq. (54) multiplied by $u^i$:

$$\frac{\partial}{\partial t} \left( \sqrt{g} H u^i \right) + \frac{\partial}{\partial x^j} \left( \sqrt{g} H u^i u^j \right) = -2\sqrt{g} \, \Gamma^i_{j0} H u^j - \sqrt{g} \, \Gamma^i_{jk} H u^j u^k$$
$$- \frac{1}{2} \sqrt{g} \, g_r h^{ij} \frac{\partial H^2}{\partial x^j} - \sqrt{g} H g_r h^{ij} \frac{\partial h_B}{\partial x^j}, \tag{55}$$

where $i = 1, 2$ and $j, k$ are summed from 1 to 2. From the identity $\partial(\sqrt{g} h^{ij})/\partial x^j = -\sqrt{g} h^{jk} \Gamma^i_{jk}$, Eq. (55) may be rewritten as

$$\frac{\partial}{\partial t} \left( \sqrt{g} H u^i \right) + \frac{\partial}{\partial x^j} \left( \sqrt{g} \left[ H u^i u^j + \frac{1}{2} g_r h^{ij} H^2 \right] \right) = -2\sqrt{g} \, \Gamma^i_{j0} H u^j$$
$$- \sqrt{g} \, \Gamma^i_{jk} \left( H u^j u^k + \frac{1}{2} g_r h^{jk} H^2 \right) - \sqrt{g} H g_r h^{ij} \frac{\partial h_B}{\partial x^j}. \tag{56}$$

Note that Eqs. (54) and (56) may be rewritten as a single expression in a (2+1)-dimensional formalism:

$$\frac{\partial}{\partial x^\nu} (\sqrt{g} T^{\mu\nu}) = -2\sqrt{g} \, \Gamma^\mu_{j0} T^{j0} - \sqrt{g} \, \Gamma^\mu_{jk} T^{jk} - \sqrt{g} H g_r h^{\mu j} \frac{\partial h_B}{\partial x^j}, \tag{57}$$

where

$$T^{\mu\nu} = H u^\mu u^\nu + \frac{1}{2} g_r h^{\mu\nu} H^2 \tag{58}$$

is the mass-momentum tensor for the shallow-water equations. Repeated Greek indices are summed from 0 to 2.

## Appendix B. Metric terms associated with the rotated cubed-sphere grid

Consider a given panel $p$ (out of 6 identical ones). Its coordinates $x^1$ and $x^2$ are great circles intersecting at right angle at the panel's center, which is assumed to be located at geographical longitude $\lambda_p$ and geographical latitude $\phi_p$. The coordinate line $x^1 = 0$ is assumed to be rotated clockwise by an angle $\alpha_p$ with respect to a geographical meridian (oriented northward) when looking at the panel's center directly from above. Define the quantities $X = \tan x^1$, $Y = \tan x^2$, and $\delta^2 = 1 + X^2 + Y^2$. The transformation from geographical longitude and latitude $(\lambda, \phi)$ to cubed-sphere coordinates $(x^1, x^2)$ is performed via the relations

$$X = \frac{\sin(\lambda - \lambda_p)\cos\alpha_p + \sin\phi_p \cos(\lambda - \lambda_p)\sin\alpha_p - \tan\phi \cos\phi_p \sin\alpha_p}{\cos\phi_p \cos(\lambda - \lambda_p) + \tan\phi \sin\phi_p}, \tag{59}$$

$$Y = \frac{\sin(\lambda - \lambda_p)\sin\alpha_p - \sin\phi_p \cos(\lambda - \lambda_p)\cos\alpha_p + \tan\phi \cos\phi_p \cos\alpha_p}{\cos\phi_p \cos(\lambda - \lambda_p) + \tan\phi \sin\phi_p}. \tag{60}$$

Each panel of the cubed-sphere grid is a gnomonic projection from the 2-sphere to the plane tangent to the panel's center. For instance, a panel with $\lambda_p = 0$ and $\alpha_p = 0$ centered at the north pole ($\phi_p = \pi/2$) has $\tan x^1 = y/z$ and $\tan x^2 = -x/z$, where $(x, y, z)$ are 3-dimensional Cartesian coordinates of the Euclidean space in which the 2-sphere is embedded.

Tensors and Christoffel symbols are transformed following the usual laws (see for instance Eqs. (6), (7) and (17) in [75]). It may be shown that the covariant metric tensor in 2+1 dimensions is

$$g_{00} = 1 + \frac{a^2}{\delta^2}\Omega^2 \left(\delta^2 - [\sin\phi_p - X\cos\phi_p \sin\alpha_p + Y\cos\phi_p \cos\alpha_p]^2\right), \tag{61}$$

$$g_{01} = \frac{a^2}{\delta^2}\Omega(1+X^2)\left(\cos\phi_p \cos\alpha_p - Y\sin\phi_p\right) = g_{10}, \tag{62}$$

$$g_{02} = \frac{a^2}{\delta^2}\Omega(1+Y^2)\left(\cos\phi_p \sin\alpha_p + X\sin\phi_p\right) = g_{20}, \tag{63}$$

$$g_{11} = \frac{a^2}{\delta^4}(1+X^2)^2(1+Y^2), \tag{64}$$

$$g_{12} = -\frac{a^2}{\delta^4}XY(1+X^2)(1+Y^2) = g_{21}, \tag{65}$$

$$g_{22} = \frac{a^2}{\delta^4}(1+X^2)(1+Y^2)^2. \tag{66}$$

The spatial metric tensor $h^{ij}$, where $h^{ij}g_{jk} = \delta^i_k$, takes the form

$$h^{11} = \frac{\delta^2}{a^2(1+X^2)}, \tag{67}$$

$$h^{12} = \frac{XY\delta^2}{a^2(1+X^2)(1+Y^2)} = h^{21}, \tag{68}$$

$$h^{22} = \frac{\delta^2}{a^2(1+Y^2)}. \tag{69}$$

The term $\sqrt{g}$ is most simply calculated as the inverse of the square root of the determinant of $h^{ij}$:

$$\sqrt{g} = \frac{a^2(1+X^2)(1+Y^2)}{\delta^3}. \tag{70}$$

The Christoffel symbols of the second kind with mixed space-time components depend on the position of the panel because the rotation at constant angular velocity $\Omega$ around a given axis breaks the (space-time) symmetry. These Christoffel symbols associated with the Coriolis acceleration may be calculated based on a space-time tensor formalism, as in [75]. They are

$$\Gamma^1_{01} = \frac{\Omega XY}{\delta^2}\left(\sin\phi_p - X\cos\phi_p \sin\alpha_p + Y\cos\phi_p \cos\alpha_p\right) = \Gamma^1_{10}, \tag{71}$$

$$\Gamma^1_{02} = -\frac{\Omega(1+Y^2)}{\delta^2}(\sin\phi_p - X\cos\phi_p \sin\alpha_p + Y\cos\phi_p \cos\alpha_p) = \Gamma^1_{20}, \tag{72}$$

$$\Gamma^2_{01} = \frac{\Omega(1+X^2)}{\delta^2}(\sin\phi_p - X\cos\phi_p \sin\alpha_p + Y\cos\phi_p \cos\alpha_p) = \Gamma^2_{10}, \tag{73}$$

$$\Gamma^2_{02} = -\frac{\Omega XY}{\delta^2}\left(\sin\phi_p - X\cos\phi_p \sin\alpha_p + Y\cos\phi_p \cos\alpha_p\right) = \Gamma^2_{20}. \tag{74}$$

The spatial components of the Christoffel symbols of the second kind are

$$\Gamma^1_{11} = \frac{2XY^2}{\delta^2}, \tag{75}$$

$$\Gamma^1_{12} = -\frac{Y(1+Y^2)}{\delta^2} = \Gamma^1_{21}, \tag{76}$$

$$\Gamma^1_{22} = 0, \tag{77}$$

$$\Gamma^2_{11} = 0, \tag{78}$$

$$\Gamma^2_{12} = -\frac{X(1+X^2)}{\delta^2} = \Gamma^2_{21}, \tag{79}$$

$$\Gamma^2_{22} = \frac{2X^2 Y}{\delta^2}. \tag{80}$$

Because of spherical symmetry, the tensor $h^{\mu\nu}$, the quantity $\sqrt{g}$ and the spatial components of the Christoffel symbols of the second kind have the same form on all six panels.

The 2-dimensional cubed-sphere coordinates have an interesting property: it may be shown that the contraction $h^{jk}\Gamma^i_{jk}$ ($i=1,2$ and $j,k$ are summed from 1 to 2) vanishes (see also [25]), thus allowing a simplification of Eq. (2). Note however that $h^{jk}\Gamma^i_{jk}$ does not vanish in general and that this term cannot be discarded from the equations of motion when other coordinates are employed.

A general rotation of the global cubed-sphere grid may be characterized by three angles which are here assumed to be $(\lambda_0, \phi_0, \alpha_0)$ associated with panel 0. The other angles $(\lambda_p, \phi_p, \alpha_p)$ associated with panels 1 to 5 are obtained as functions of these $(\lambda_0, \phi_0, \alpha_0)$. They are found to be

$$\lambda_1 = \tan^{-1}\left(\frac{\cos\lambda_0 \cos\alpha_0 + \sin\lambda_0 \sin\phi_0 \sin\alpha_0}{\cos\lambda_0 \sin\phi_0 \sin\alpha_0 - \sin\lambda_0 \cos\alpha_0}\right), \tag{81}$$

$$\phi_1 = -\sin^{-1}(\cos\phi_0 \sin\alpha_0), \tag{82}$$

$$\alpha_1 = \tan^{-1}\left(\frac{\sin\phi_0}{\cos\phi_0 \cos\alpha_0}\right), \tag{83}$$

$$\lambda_2 = \lambda_0 + \pi, \tag{84}$$

$$\phi_2 = -\phi_0, \tag{85}$$

$$\alpha_2 = -\alpha_0, \tag{86}$$

$$\lambda_3 = -\tan^{-1}\left(\frac{\cos\lambda_0 \cos\alpha_0 + \sin\lambda_0 \sin\phi_0 \sin\alpha_0}{\sin\lambda_0 \cos\alpha_0 - \cos\lambda_0 \sin\phi_0 \sin\alpha_0}\right), \tag{87}$$

$$\phi_3 = \sin^{-1}(\cos\phi_0 \sin\alpha_0), \tag{88}$$

$$\alpha_3 = -\tan^{-1}\left(\frac{\sin\phi_0}{\cos\phi_0 \cos\alpha_0}\right), \tag{89}$$

$$\lambda_4 = \tan^{-1}\left(\frac{\cos\lambda_0 \sin\alpha_0 - \sin\lambda_0 \sin\phi_0 \cos\alpha_0}{-\cos\lambda_0 \sin\phi_0 \cos\alpha_0 - \sin\lambda_0 \sin\alpha_0}\right), \tag{90}$$

$$\phi_4 = \sin^{-1}(\cos\phi_0 \cos\alpha_0), \tag{91}$$

$$\alpha_4 = \tan^{-1}\left(\frac{\cos\phi_0 \sin\alpha_0}{-\sin\phi_0}\right), \tag{92}$$

$$\lambda_5 = \tan^{-1}\left(\frac{\sin\lambda_0 \sin\phi_0 \cos\alpha_0 - \cos\lambda_0 \sin\alpha_0}{\cos\lambda_0 \sin\phi_0 \cos\alpha_0 + \sin\lambda_0 \sin\alpha_0}\right), \tag{93}$$

$$\phi_5 = -\sin^{-1}(\cos\phi_0 \cos\alpha_0), \tag{94}$$

$$\alpha_5 = \tan^{-1}\left(\frac{\cos\phi_0 \sin\alpha_0}{\sin\phi_0}\right). \tag{95}$$

Care must be taken when choosing the correct branch of the arctangent operators for various panels. In the special case where $\phi_0 = \alpha_0 = 0$, different formulas for $\lambda_4, \alpha_4, \lambda_5$ and $\alpha_5$ must be used due to the singularity at the poles. In this case, $\lambda_4 = 0 = \lambda_5$ and $\alpha_4 = -\lambda_0 = -\alpha_5$. Examples are provided in the accompanying implementation code (see §7).

## Appendix C. Consistency relations at the interfaces of panels

At the interface of two panels, consistency relations on tensor components may be established. Consider for instance the contravariant vector components $A^1_{(0)}$, $A^2_{(0)}$ on panel 0 and their relations to the contravariant vector components $A^1_{(1)}$, $A^2_{(1)}$ on panel 1. Note that due to spherical symmetry, the form of these relations is not altered by rotating the grid. One may then consider the case $\lambda_0 = \phi_0 = \alpha_0 = \phi_1 = \alpha_1 = 0$ and $\lambda_1 = \pi/2$. From $X_{(0)} = \tan\lambda$, $Y_{(0)} = \tan\phi/\cos\lambda$, $X_{(1)} = -1/\tan\lambda$ and $Y_{(1)} = \tan\phi/\sin\lambda$, the transformation laws from coordinates on panel 0 to coordinates on panel 1 lead to

$$\frac{\partial x^1_{(0)}}{\partial x^1_{(1)}} = 1, \tag{96}$$

$$\frac{\partial x^1_{(0)}}{\partial x^2_{(1)}} = 0, \tag{97}$$

$$\frac{\partial x^2_{(0)}}{\partial x^1_{(1)}} = \frac{Y_{(1)}(1 + X^2_{(1)})}{X^2_{(1)}(1 + Y^2_{(1)})} \quad \Rightarrow \quad \frac{\partial x^2_{(0)}}{\partial x^1_{(1)}} = \frac{2Y_{(1)}}{1 + Y^2_{(1)}} \text{ at } X_{(1)} = -1, \tag{98}$$

$$\frac{\partial x^2_{(0)}}{\partial x^2_{(1)}} = -\frac{X_{(1)}(1 + Y^2_{(1)})}{X^2_{(1)} + Y^2_{(1)}} \quad \Rightarrow \quad \frac{\partial x^2_{(0)}}{\partial x^2_{(1)}} = 1 \text{ at } X_{(1)} = -1. \tag{99}$$

These relations are used at the interface $X_{(1)} = -1$ to convert contravariant vector components from panel 1 to panel 0. A similar approach is used to convert covariant vector components. This procedure may be performed at the 12 interfaces of the cubed-sphere grid.

The consistency relations for contravariant vector components $A^i_{(p)}$ and covariant vector components $A_{(p)i}$ on panel $p$ at the 12 interfaces are provided below. At the interface of panels $(p, q) = (0, 1), (1, 2), (2, 3), (3, 0)$, one obtains

$$A^1_{(p)} = A^1_{(q)}, \tag{100}$$

$$A^2_{(p)} = \frac{2Y}{1 + Y^2} A^1_{(q)} + A^2_{(q)}, \tag{101}$$

$$A_{(p)1} = A_{(q)1} - \frac{2Y}{1 + Y^2} A_{(q)2}, \tag{102}$$

$$A_{(p)2} = A_{(q)2}. \tag{103}$$

At the interface of panels $(4, 0)$, one obtains

$$A^1_{(4)} = A^1_{(0)} - \frac{2X}{1 + X^2} A^2_{(0)}, \tag{104}$$

$$A^2_{(4)} = A^2_{(0)}, \tag{105}$$

$$A_{(4)1} = A_{(0)1}, \tag{106}$$

$$A_{(4)2} = \frac{2X}{1 + X^2} A_{(0)1} + A_{(0)2}, \tag{107}$$

with $X$ defined on panel 0. At the interface of panels $(4, 1)$, one obtains

$$A^1_{(4)} = -A^2_{(1)}, \tag{108}$$

$$A^2_{(4)} = A^1_{(1)} - \frac{2X}{1 + X^2} A^2_{(1)}, \tag{109}$$

$$A_{(4)1} = -\frac{2X}{1 + X^2} A_{(1)1} - A_{(1)2}, \tag{110}$$

$$A_{(4)2} = A_{(1)1}, \tag{111}$$

with $X$ defined on panel 1. At the interface of panels $(4, 2)$, one obtains

$$A^1_{(4)} = -A^1_{(2)} + \frac{2X}{1 + X^2} A^2_{(2)}, \tag{112}$$

$$A^2_{(4)} = -A^2_{(2)}, \tag{113}$$

$$A_{(4)1} = -A_{(2)1}, \tag{114}$$

$$A_{(4)2} = -\frac{2X}{1 + X^2} A_{(2)1} - A_{(2)2}, \tag{115}$$

with $X$ defined on panel 2. At the interface of panels $(4, 3)$, one obtains

$$A^1_{(4)} = A^2_{(3)}, \tag{116}$$

$$A^2_{(4)} = -A^1_{(3)} + \frac{2X}{1 + X^2} A^2_{(3)}, \tag{117}$$

$$A_{(4)1} = \frac{2X}{1+X^2} A_{(3)1} + A_{(3)2}, \tag{118}$$

$$A_{(4)2} = -A_{(3)1}, \tag{119}$$

with $X$ defined on panel 3. At the interface of panels $(5, 0)$, one obtains

$$A^1_{(5)} = A^1_{(0)} + \frac{2X}{1+X^2} A^2_{(0)}, \tag{120}$$

$$A^2_{(5)} = A^2_{(0)}, \tag{121}$$

$$A_{(5)1} = A_{(0)1}, \tag{122}$$

$$A_{(5)2} = -\frac{2X}{1+X^2} A_{(0)1} + A_{(0)2}, \tag{123}$$

with $X$ defined on panel 0. At the interface of panels $(5, 1)$, one obtains

$$A^1_{(5)} = A^2_{(1)}, \tag{124}$$

$$A^2_{(5)} = -A^1_{(1)} - \frac{2X}{1+X^2} A^2_{(1)}, \tag{125}$$

$$A_{(5)1} = -\frac{2X}{1+X^2} A_{(1)1} + A_{(1)2}, \tag{126}$$

$$A_{(5)2} = -A_{(1)1}, \tag{127}$$

with $X$ defined on panel 1. At the interface of panels $(5, 2)$, one obtains

$$A^1_{(5)} = -A^1_{(2)} - \frac{2X}{1+X^2} A^2_{(2)}, \tag{128}$$

$$A^2_{(5)} = -A^2_{(2)}, \tag{129}$$

$$A_{(5)1} = -A_{(2)1}, \tag{130}$$

$$A_{(5)2} = \frac{2X}{1+X^2} A_{(2)1} - A_{(2)2}, \tag{131}$$

with $X$ defined on panel 2. At the interface of panels $(5, 3)$, one obtains

$$A^1_{(5)} = -A^2_{(3)}, \tag{132}$$

$$A^2_{(5)} = A^1_{(3)} + \frac{2X}{1+X^2} A^2_{(3)}, \tag{133}$$

$$A_{(5)1} = \frac{2X}{1+X^2} A_{(3)1} - A_{(3)2}, \tag{134}$$

$$A_{(5)2} = A_{(3)1}, \tag{135}$$

with $X$ defined on panel 3.

Consistency relations for higher-rank tensors may be obtained from these rules. For instance, if a symmetric contravariant second-rank tensor is considered at a point on an interface, then $T^{00}$ remains identical on both panels; $T^{0i}$ and $T^{i0}$ transform as $A^i$; and $T^{ij}$ as the product $A^i A^j$.

### Appendix D. AUSM Riemann solver for the shallow-water equations

The Advection Upstream Splitting Method (AUSM) is a simple flux splitting method. It works by splitting the advective component of the flux from the pressure component. This appendix presents a brief overview of the elements that are relevant for its implementation in the context of the shallow-water equations. The reader is referred to [41] for more details.

From Eq. (58), define a tensor density $q^\mu \equiv \sqrt{g}\, T^{\mu 0} = \sqrt{g} H u^\mu$. One may express $\sqrt{g}\, T^{\mu\nu}$ as a function of $q^\mu$:

$$\sqrt{g}\, T^{\mu\nu} = \frac{q^\mu q^\nu}{q^0} + \frac{1}{2\sqrt{g}} g_r h^{\mu\nu} \left(q^0\right)^2. \tag{136}$$

In AUSM, the components of the tensor densities appearing in the spatial derivative of Eq. (57) are splitted into

$$\sqrt{g}\, T^{\mu i} = M^i A^{\mu i} + P^{\mu i}, \tag{137}$$

(no sum over $i$) where $A^{\mu i} = q^\mu \sqrt{g} H a^i / q^0$ is the advective component, $M^i = u^i / a^i$ is the Froude number, $a^i$ is the intrinsic gravity wave velocity (derived in Appendix E) and $P^{\mu i} = g_r h^{\mu i} (q^0)^2 / (2\sqrt{g})$ is the so-called pressure component.

The Froude number and pressure component are decomposed into a part evaluated on one side of the interface (denoted by a "+" superscript) and another part evaluated on the other side (denoted by a "−" superscript). Specifically,

$$M^i = (M^i)^+ + (M^i)^- \tag{138}$$

and

$$P^{\mu i} = (P^{\mu i})^+ + (P^{\mu i})^-, \tag{139}$$

where $(M^i)^\pm$ and $(P^{\mu i})^\pm$ are respectively defined by Eq. (6) and Eq. (8) of [41]. Analogously to the original AUSM formulation, the Froude number is splitted a second time as follows

$$M^i = \max[0, (M^i)^+ + (M^i)^-] + \min[0, (M^i)^+ + (M^i)^-]. \tag{140}$$

Then the AUSM method with "double Froude number splitting" is

$$\sqrt{g} T^{\mu i} = \left[ (\sqrt{g} T^{\mu i})^+ + (\sqrt{g} T^{\mu i})^- \right] \tag{141}$$

where

$$(\sqrt{g} T^{\mu i})^+ = \max[0, (M^i)^+ + (M^i)^-] (A^{\mu i})^+ + (P^{\mu i})^+, \tag{142}$$

$$(\sqrt{g} T^{\mu i})^- = \min[0, (M^i)^+ + (M^i)^-] (A^{\mu i})^- + (P^{\mu i})^-. \tag{143}$$

**Appendix E. Gravity wave velocities**

The components $(\nu, \alpha)$ of the three flux Jacobian matrices $\mathcal{J}^\mu$ are defined as

$$\left(\mathcal{J}^\mu\right)^\nu_\alpha \equiv \frac{\partial(\sqrt{g} T^{\mu\nu})}{\partial q^\alpha}. \tag{144}$$

Note that the left-hand side of Eq. (57) may be rewritten with the flux Jacobian matrices as

$$\frac{\partial}{\partial x^\nu}(\sqrt{g} T^{\mu\nu}) = \left(\mathcal{J}^\mu\right)^\nu_\alpha \frac{\partial q^\alpha}{\partial x^\nu}. \tag{145}$$

They are written explicitly as $\mathcal{J}^0 = I_{3\times 3}$ (the identity matrix) and

$$\mathcal{J}^1 = \begin{pmatrix} 0 & 1 & 0 \\ g_r h^{11} H - u^1 u^1 & 2u^1 & 0 \\ g_r h^{12} H - u^1 u^2 & u^2 & u^1 \end{pmatrix}; \quad \mathcal{J}^2 = \begin{pmatrix} 0 & 0 & 1 \\ g_r h^{12} H - u^1 u^2 & u^2 & u^1 \\ g_r h^{22} H - u^2 u^2 & 0 & 2u^2 \end{pmatrix}. \tag{146}$$

The three eigenvalues of $\mathcal{J}^0$ correspond to $u^0 = 1$. The three eigenvalues of $\mathcal{J}^1$ are

$$u^1, \; u^1 \pm \sqrt{h^{11} g_r H}, \tag{147}$$

and those of $\mathcal{J}^2$ are

$$u^2, \; u^2 \pm \sqrt{h^{22} g_r H}. \tag{148}$$

Since the shallow-water equations are hyperbolic, the eigenvalues (147) and (148) are real and distinct. The intrinsic gravity wave velocities are deduced from these eigenvalues as $a^1 \equiv \sqrt{h^{11} g_r H}$ and $a^2 \equiv \sqrt{h^{22} g_r H}$.

S. Gaudreault, M. Charron, V. Dallerit et al.